\documentclass[11pt]{article}
\usepackage[T1]{fontenc}
\usepackage{lmodern,amsmath,amssymb,amsthm,mathtools,booktabs,array,enumitem}
\usepackage[margin=1in]{geometry}
\usepackage{natbib}
\usepackage{xcolor,graphicx,hyperref,microtype}
\definecolor{linkblue}{HTML}{183C62}
\hypersetup{colorlinks=true,linkcolor=linkblue,citecolor=linkblue,urlcolor=linkblue,pdftitle={Exact Last-Iterate Bounds and Optimal Stepsizes for Proximal Stochastic Gradient Descent},pdfauthor={Laurent Condat and Peter Richtarik}}
\newtheorem{theorem}{Theorem}[section]
\newtheorem{lemma}[theorem]{Lemma}
\newtheorem{proposition}[theorem]{Proposition}
\newtheorem{corollary}[theorem]{Corollary}
\theoremstyle{definition}
\newtheorem{assumption}[theorem]{Assumption}
\newtheorem{remark}[theorem]{Remark}
\newcommand{\R}{\mathbb R}
\newcommand{\E}{\mathbb E}
\newcommand{\prox}{\operatorname{prox}}
\newcommand{\dom}{\operatorname{dom}}
\newcommand{\ip}[2]{\langle #1,#2\rangle}
\newcommand{\norm}[1]{\lVert #1\rVert}
\newcommand{\gap}{\Delta_T}
\newcommand{\cF}{\mathcal F}
\newcommand{\cW}{\mathcal W}
\newcommand{\Ec}[1]{\E\!\left[#1\mid\cF^t\right]}
\newcommand{\eqdef}{\coloneqq}
\newcommand{\xc}{\mathring{x}}
\newcommand{\qc}{\mathring{q}}
\newcommand{\gc}{\check{\gamma}}
\newcommand{\hub}{\operatorname{hub}}
\DeclareMathOperator*{\argmin}{arg\,min}
\newcommand{\cH}{\mathcal{H}}
\definecolor{mydarkblue}{rgb}{0,0.0,0.8}
\newcommand{\ff}{\mathcolor{mydarkblue}{f}}
\definecolor{mydarkgreen}{rgb}{0,0.45,0}
\newcommand\RR{\mathcolor{mydarkgreen}{R}}
\usepackage{xspace}
\usepackage{scalefnt}
\definecolor{mydarkredd}{rgb}{0.45,0.16,0.0}
\newcommand{\algn}[1]{{\sf\color{mydarkredd}\scalefont{0.97}{#1}}\xspace}
\newcommand{\algno}{\algn{Prox-SGD}}

\allowdisplaybreaks[1]

\setlist{itemsep=3pt,topsep=4pt}
\newcolumntype{P}[1]{>{\raggedright\arraybackslash}p{#1}}

\title{\textbf{Exact Last-Iterate Bounds and Optimal Stepsizes\\ for Proximal Stochastic Gradient Descent}}
\author{Laurent Condat\ \&\  Peter Richt\'arik\\[2pt] King Abdullah University of Science and Technology (KAUST)}
\date{November 23, 2026}

\begin{document}
\maketitle
\begin{abstract}
Proximal Stochastic Gradient Descent 
 minimizes the sum of an $L$-smooth convex function $f$, known through unbiased stochastic gradients of bounded variance, and a nonsmooth convex function $R$, known through its proximity operator. We derive the exact worst-case value of the expected objective gap of its last iterate, for every predetermined stepsize schedule in a broad class, which contains all the schedules with stepsizes at most $\frac{3}{2L}$. 
This extends to the composite setting a recent result of Grimmer and Shen for the smooth case $R=0$. 
The value is the sum of an initial-error term and of one term per iteration, which is the price of the noise injected at this iteration, damped by the stepsizes that remain after it. The upper bound is proved with a simple potential, the Moreau envelope of the objective with a smoothing parameter proportional to the remaining stepsize sum, and it is attained by an explicit instance. The formula shows that a logarithmic factor is unavoidable with a constant stepsize,
and it gives the exact price of the proximity operator, relative to the smooth case.
It also makes the choice of the stepsizes a tractable optimization problem: for a uniform variance bound $\sigma^2$, the best schedule is unique, decreasing, and computed by a scalar search, and it is minimax optimal among all predetermined schedules, unless the noise level is low. Asymptotically, the minimax value $M_T$ after $T$ iterations satisfies $(D\sigma/M_T)^2=T+\log T+C_\star+o(1)$, where $D$ is a bound on the initial distance to a solution and $C_\star\approx1.768$ is a universal constant; unlike in the smooth case, $L$ plays no role at this order. This law is achieved by an explicit schedule, which is the optimal one when $f=0$. 
\end{abstract}
\noindent\textbf{Keywords:} stochastic optimization; proximal gradient; last iterate; exact worst-case analysis; stepsize design.

\tableofcontents

\section{Introduction}\label{sec:intro}
Many problems in signal processing, statistics, and machine learning consist in minimizing the sum of a smooth data-fidelity term and a nonsmooth regularizer or constraint. Proximal gradient methods, which alternate a gradient step on the smooth term and a proximal step on the nonsmooth one, are the basic tools for such problems. When the gradient of the smooth term is only available through a stochastic estimate, the resulting algorithm is \emph{Proximal Stochastic Gradient Descent} (\algno), a method so simple and so widely used that one could expect its behavior to be completely understood. 
It is not: the exact convergence speed of its last iterate, the point actually returned after a given number of iterations, was not known, even for a constant stepsize. In this paper, we derive new convergence rates and optimal stepsize schedules for \algno with stochastic gradient errors of bounded variance.
Our setting is the convex optimization problem:
\begin{equation*}
\mathrm{Find}\ \,x^\star \in \argmin_{x\in\cH}\, \Psi(x)\eqdef \ff(x)+\RR(x),
\end{equation*}
under the standing assumption:
\begin{assumption}\label{ass:model}
$\cH$ is a real Hilbert space. 
The function $\ff:\cH\to\R$ is convex and differentiable, with $L$-Lipschitz gradient, for some $L>0$. The function $\RR:\cH\to \R\cup\{+\infty\}$ is proper, lower semicontinuous, and convex \citep{bau17}. 
We let $x^\star\in\arg\min\Psi$, a set that is assumed to be nonempty.
The initial point $x^0\in\cH$ is deterministic and fixed ($x^0\notin \dom \RR$ is permitted), and $\norm{x^0-x^\star}\leq D$ for some $D\geq0$. We write $\Psi^\star\eqdef\Psi(x^\star)$.
\end{assumption}
In words, $\ff$ is $L$-smooth and
$\RR$ is a possibly nonsmooth function, which is only required to be \emph{proximable}:
its proximity operator (prox, for short) is simple and can be evaluated exactly, as is the case for many functions of practical interest \citep{par14};
we recall that the prox of $\gamma \RR$, for any $\gamma>0$, is the mapping $\mathrm{prox}_{\gamma \RR}:x\in\cH \mapsto \arg\min_{y\in\cH}\big\{\gamma\RR(y)+\frac1{2}\|y-x\|^2\big\}$.

Let the number $T\ge 1$ of iterations be fixed. Given
a sequence $(\gamma_t)_{t=0}^{T-1}$
of positive stepsizes and the initial point $x^0$,
\algno 
 iterates \citep{ros20}
\begin{equation}
x^{t+1} \eqdef \mathrm{prox}_{\gamma_t \RR} \Big(x^t-\gamma_t \big(\nabla \ff (x^t)+e^t\big) \Big) ,\quad \forall t=0,\ldots, T-1,\label{eq:algorithm}
\end{equation}
where the error $e^t\in \cH$ is a square-integrable random vector, 
which is conditionally unbiased, with bounded variance; that is, denoting by $\cF^t$ the sigma-algebra generated by the iterates $x^0,\ldots,x^t$, and by $\sigma_t^2\ge 0$ a prescribed variance bound,
\begin{equation}
\Ec{e^t}=0,\quad  \E\big[\|e^t\|^2\big]\le\sigma_t^2,\quad \forall t=0,\ldots,T-1.\label{eq:noise}
\end{equation}
Neither the independence of the errors $e^0,\ldots,e^{T-1}$ nor an almost-sure bound on the conditional variance $\Ec{\|e^t\|^2}$ is required:
only the conditional unbiasedness and the unconditional second-moment bounds in \eqref{eq:noise} enter the proofs. Conversely, the worst-case instances constructed in Section~\ref{sec:main} satisfy $\Ec{\|e^t\|^2}=\sigma_t^2$ almost surely, so that the exact values obtained below are the same under this stronger conditional model.
We define the objective gap
\begin{equation*}
\Delta_t\eqdef \Psi(x^t)-\Psi^\star\ \ge 0,\quad \forall t=0,\ldots,T.
\end{equation*}
Our goal is to characterize the last-iterate expected gap $\E[\Delta_T]$ of \algno,
and to find stepsize schedules $\boldsymbol{\gamma}\eqdef(\gamma_t)_{t=0}^{T-1}$ that minimize this gap, by deriving exact bounds, attained by worst-case instances, over the class of problems satisfying Assumption~\ref{ass:model}, for given $D$, $L$, $T$, and $\boldsymbol{\sigma}\eqdef(\sigma_t)_{t=0}^{T-1}$.
The schedule is \emph{predetermined}: the stepsizes are deterministic, chosen before the run, possibly as functions of the horizon $T$. A zero stepsize, understood as $x^{t+1}=x^t$, need not be considered: such an iteration can be removed, with the horizon reduced accordingly.\medskip

\noindent\textbf{The last-iterate question.\ } As recalled in Section~\ref{sec:classical}, standard Lyapunov arguments,
with a constant stepsize $\gamma\le\frac1L$ and a uniform variance bound $\sigma^2$, yield the bound $\frac{D^2}{2\gamma T}+\gamma\sigma^2$ on the expected gap of the \emph{averaged} iterate, and nothing on the last iterate, since the objective gaps are not monotone in the presence of noise.
Yet the last iterate is what the algorithm returns, and averaging is unwelcome in the composite setting, where it destroys the structure, such as sparsity, enforced by the prox.
Under the noise model \eqref{eq:noise}, last-iterate bounds for \algno have been obtained only recently \citep{lz24};
with a constant stepsize, their noise term is of order $\gamma\sigma^2\log T$ instead of $\gamma\sigma^2$: the noise injected at a given iteration is felt by the last iterate with a weight decaying like the inverse of the number of remaining iterations, and the sum of these weights is a harmonic number. In the smooth case $\RR=0$, \citet{gs26} recently established the exact worst-case value of $\E[\Delta_T]$ for Stochastic Gradient Descent (\algn{SGD}), for every schedule of stepsizes in $[0,\frac{3}{2L}]$; this shows in particular that the logarithm is unavoidable with a constant stepsize, and the exact formula allowed them to derive the minimax optimal schedule, which is decreasing and removes the logarithm.\medskip

\noindent\textbf{Our results in words.\ }
We determine exactly how large the expected gap of the last iterate can be, for every schedule in a broad class, which contains all the schedules with stepsizes at most $\frac{3}{2L}$. With $S_t$ denoting the sum of the stepsizes $\gamma_t,\ldots,\gamma_{T-1}$ that remain at iteration $t$, the answer is
\begin{equation}
\E[\Delta_T]\ \le\ \underbrace{\frac{D^2}{4S_0}}_{\text{initial distance}}\ +\ \sum_{t=0}^{T-1}\ \underbrace{\frac{(\gamma_t\sigma_t)^2}{2d_t}}_{\text{noise of iteration }t}\,,\qquad\text{where}\quad d_t\eqdef 2S_{t+1}+\omega(\gamma_t)
\label{eq:inwords}
\end{equation}
and $\omega(\gamma_t)$ is a positive number, at most $2\gamma_t$, defined in \eqref{eq:omega} below. Each term is a price. The first one is the price of starting at distance $D$ from a solution; it is the exact rate of the method without noise, and it decreases when the total stepsize $S_0$ increases. The term of index $t$ is the price of the noise injected at iteration $t$: a perturbation of magnitude $\gamma_t\sigma_t$, damped by the stepsizes that remain after it, through $S_{t+1}$, and by the step itself, through $\omega(\gamma_t)$. Nothing is lost in this bound: for every such schedule, there is a problem, in dimension at least $T+1$, for which \eqref{eq:inwords} is an equality.

Three conclusions follow. First, with a constant stepsize, the last injections are hardly damped, and the prices add up to a logarithm of $T$, which is therefore unavoidable. Second, a decreasing schedule, with small steps when little damping remains, avoids this logarithm; since the formula is exact,
the best schedule is the one that minimizes the right-hand side of \eqref{eq:inwords}:
it is unique, it is decreasing, it is computed by a scalar search, and no predetermined schedule has a better guarantee, unless the noise level is low. Third, for a long horizon, this best guarantee is $D\sigma/\sqrt{T+\log T+C_\star+o(1)}$, with a universal constant $C_\star\approx1.768$;
it does not depend on $L$, unlike that of \algn{SGD} at this order,
and it is achieved by a simple explicit schedule, which depends on $T$ and on the ratio $D/\sigma$ only.\medskip

\noindent\textbf{The envelope idea and the composite obstacle.\ }
The key to our analysis is a potential adapted to the horizon: the Moreau envelope of the objective at the current iterate, with a smoothing parameter proportional to the remaining stepsize sum, which vanishes at the horizon, where the envelope is the objective itself. The noise injected at iteration $t$ is felt by this potential through its curvature, which grows from nearly $0$, when the horizon is far, to its maximum value at the end, and the sum of these curvatures over the injections is the logarithm. In the smooth case $\RR=0$, the envelope appears in \citet[Section~2.1]{gs26} as an interpretation of their inductive argument, which is formulated with an auxiliary sequence of points.
When $\RR\neq0$, the noise is inside the prox, which does not average it out, and the envelope of $\Psi$ loses its improved curvature: the argument breaks down, at first order in the noise, as shown by a one-dimensional example in Section~\ref{sec:noisystep}.
\emph{Our analysis makes the envelope itself the potential}, without any auxiliary sequence, and controls the whole noisy proximal step by a single pointwise inequality, Lemma~\ref{lem:local}: the envelope after the step is bounded by the envelope before the step, plus a term linear in the noise with a predictable coefficient, which vanishes in conditional expectation, plus a quadratic term whose coefficient is the exact price of the injection. The inequality is proved by an explicit sum-of-squares identity (Appendix~\ref{app:certificate}), and its noise coefficient cannot be improved (Proposition~\ref{prop:localsharp}).\medskip

\noindent\textbf{Notation.\ }
For a schedule $\boldsymbol{\gamma}$, define the remaining stepsize sums
\begin{equation}
 S_t\eqdef\sum_{j=t}^{T-1}\gamma_j,\quad \forall t=0,\ldots,T-1,\qquad S_T\eqdef0.\label{eq:tails}
\end{equation}
The schedule is said to be \emph{admissible} if
\begin{equation}
 L\gamma_t+\frac{\gamma_t}{2S_t}\leq2,\quad \forall t=0,\ldots,T-1;\label{eq:admissible}
\end{equation}
since $\gamma_t/(2S_t)\le\frac12$, this is the case if every $\gamma_t\le\frac{3}{2L}$. 
This is the class of schedules for which \eqref{eq:inwords} is proved. The function $\omega$ that appears in it is
\begin{equation}
 \omega(\gamma)\eqdef\tfrac1L\Big(1-(1-L\gamma)_+^2\Big),\qquad\text{where}\quad(\,\cdot\,)_+\eqdef\max\{\,\cdot\,,0\},\label{eq:omega}
\end{equation}
that is, $\omega(\gamma)=\gamma(2-L\gamma)$ if $L\gamma\le1$ and $\omega(\gamma)=\frac1L$ if $L\gamma\ge1$. This function is positive, nondecreasing, and continuously differentiable on $(0,\infty)$; the role of the threshold $\frac1L$ is explained in Section~\ref{sec:noisystep}. Thus, for every positive schedule, the positive numbers $d_t$ of \eqref{eq:inwords} are
\begin{equation}
 d_t\eqdef2S_{t+1}+\omega(\gamma_t)\overset{\eqref{eq:tails},\eqref{eq:omega}}{=}\begin{cases}2S_t-L\gamma_t^2,&\text{if }L\gamma_t\le1,\\[1mm] 2S_{t+1}+\frac1L,&\text{if }L\gamma_t\ge1,\end{cases}\qquad \forall t=0,\ldots,T-1,\label{eq:denominators}
\end{equation}
and the right-hand side of \eqref{eq:inwords} is called the \emph{cost} of the schedule:
\begin{equation}
 F_{L,T}(\boldsymbol{\gamma};D,\boldsymbol{\sigma})\eqdef\frac{D^2}{4S_0}
       +\frac12\sum_{t=0}^{T-1}\frac{\sigma_t^2\gamma_t^2}{d_t}.\label{eq:cost}
\end{equation}
Let $\cW_{L,T}(\boldsymbol{\gamma};D,\boldsymbol{\sigma})$ denote the worst-case value of $\E[\Delta_T]$, that is, its supremum over all instances in $\cH$: all functions $\ff$ and $\RR$, initial points $x^0$, and noise distributions satisfying Assumption~\ref{ass:model} and \eqref{eq:noise}.
In this notation, as in those introduced later,
the subscripts record the smoothness constant $L$, which will be varied, and the horizon $T$; 
the argument before the semicolon is the schedule, the design variable; the arguments after it are the class parameters over which the worst case is taken.\medskip

\noindent\textbf{Main result.\ } Our main result is the following theorem, which makes \eqref{eq:inwords} precise. It is proved in Section~\ref{sec:main}.

\begin{theorem}[Exact performance of a prescribed schedule]\label{thm:exact}
Let $\boldsymbol{\gamma}$ be a positive predetermined schedule.
\begin{enumerate}[label=(\roman*)]
\item If $\boldsymbol{\gamma}$ is admissible, that is, satisfies \eqref{eq:admissible}, which is the case if every $\gamma_t\le\frac{3}{2L}$, then
\begin{equation}
 \E[\gap]\leq F_{L,T}(\boldsymbol{\gamma};D,\boldsymbol{\sigma}).\label{eq:upper}
\end{equation}
\item If $\dim\cH\ge T+1$, then
\begin{equation}
 \cW_{L,T}(\boldsymbol{\gamma};D,\boldsymbol{\sigma})\geq F_{L,T}(\boldsymbol{\gamma};D,\boldsymbol{\sigma}):\label{eq:universallower}
\end{equation}
there is an instance, which can be chosen with a unique minimizer, on which $\gap=F_{L,T}(\boldsymbol{\gamma};D,\boldsymbol{\sigma})$ for every realization of its noise signs.
\end{enumerate}
Consequently, if $\boldsymbol{\gamma}$ is admissible and $\dim\cH\ge T+1$, then
\begin{equation}
 \cW_{L,T}(\boldsymbol{\gamma};D,\boldsymbol{\sigma})\overset{\eqref{eq:upper},\eqref{eq:universallower}}{=}F_{L,T}(\boldsymbol{\gamma};D,\boldsymbol{\sigma}).\label{eq:exact}
\end{equation}
\end{theorem}

\noindent\textbf{Contributions.\ } Our contributions are the following.\smallskip

\noindent 1) \emph{Exact performance of a schedule} (Section~\ref{sec:main}). Theorem~\ref{thm:exact} 
shows that the cost is an upper bound on $\E[\Delta_T]$ for every admissible schedule, and a lower bound on its worst-case value for every positive schedule; hence it is the exact worst-case value for every admissible schedule. The attaining instance is separable, made of Huber and $\ell_1$ functions of $T+1$ orthonormal coordinates, with Rademacher noise. When $\ff=0$, the cost with $L=0$ is the exact worst case for every positive schedule, without restriction on the stepsizes (Corollary~\ref{cor:ppa}).\smallskip

\noindent 2) \emph{The price of the prox.} In the smooth case $\RR=0$, the exact worst case of \algn{SGD} \citep{gs26} has the same form as the cost \eqref{eq:cost}, with $4S_0+\frac2L$ in place of $4S_0$ and $\frac1L$ in place of $\omega(\gamma_t)$; see Table~\ref{tab:price}. Since $\omega(\gamma_t)=\frac1L-\frac1L(1-L\gamma_t)_+^2$, the composite case is never better: the constant $\frac2L$ is lost in the first term,
and the noise is less damped at the steps smaller than $\frac1L$. The reason is that the prox is nonlinear:
a centered noise at its input displaces the mean of its output. By exactness, this price is real; it is given by the identity \eqref{eq:price}. The first term of the cost is the exact worst case of the deterministic proximal gradient method, previously known for stepsizes at most $\frac1L$ only.\smallskip

\begin{table}[t]\centering
\caption{The exact worst-case value of $\E[\Delta_T]$ is the sum of an initial-error term and of one noise term per iteration, in the smooth case $\RR=0$ \citep{gs26} and in the composite case (this paper). The last row is the minimax value $M$ over all predetermined schedules, for a uniform variance bound $\sigma^2$ and $T\to\infty$ (under a sufficient-noise condition in the smooth case).}\label{tab:price}
\smallskip
\small
\begin{tabular}{lcc}
\toprule
 & \algn{SGD} ($\RR=0$) & \algno\\
\midrule
Initial-error term & $\dfrac{D^2}{4S_0+2/L}$ & $\dfrac{D^2}{4S_0}$\\[4mm]
Noise term of iteration $t$ & $\dfrac{\sigma_t^2\gamma_t^2}{2\big(2S_{t+1}+1/L\big)}$ & $\dfrac{\sigma_t^2\gamma_t^2}{2\big(2S_{t+1}+\omega(\gamma_t)\big)}$\\[4mm]
Exact for & \multicolumn{2}{c}{admissible schedules \eqref{eq:admissible}}\\[1mm]
Minimax law $(D\sigma/M)^2$ & $T+2\sqrt{\sigma/(LD)}\,T^{3/4}+o(T^{3/4})$ & $T+\log T+C_\star+o(1)$\\
\bottomrule
\end{tabular}
\end{table}

\noindent 3) \emph{Optimal schedules} (Section~\ref{sec:design}). For a uniform variance bound and $D,\sigma>0$, the cost is a strictly convex function over all positive schedules; its unique minimizer is strictly decreasing and is computed by a scalar search (Theorems~\ref{thm:convex}--\ref{thm:scalarconstruction}). Whenever this minimizer is admissible, which is the case if $\sigma\ge\frac{LD}{2\sqrt3}$ or if $T+3\ge\frac{4L^2D^2}{\sigma^2}$, it is the unique minimax optimal schedule among \emph{all} positive predetermined schedules (Theorem~\ref{thm:global} and Corollary~\ref{cor:globalT}): the finite-horizon design problem is then solved.
In the case $\ff=0$, the optimal schedule is given in closed form by a Bellman recursion, which also handles a nonuniform variance profile, for which stopping before the horizon can be optimal, and yields an exact decomposition of the excess cost of a suboptimal schedule (Theorems~\ref{thm:ppaopt} and~\ref{thm:heterogeneous}, Corollary~\ref{cor:schedulegap}).\smallskip

\noindent 4) \emph{Asymptotic minimax law} (Section~\ref{sec:asymptotic}).
The minimax value $M_{L,T}$ of $\E[\Delta_T]$ over \emph{all} positive predetermined schedules satisfies $(D\sigma/M_{L,T})^2=T+\log T+C_\star+o(1)$ for fixed $L,D,\sigma>0$ (Theorem~\ref{thm:secondorder}).
The optimal schedule of the case $\ff=0$, computed by a scalar recursion from $T$ and the ratio $D/\sigma$ only, achieves this law for every $L$, and linearly decaying stepsizes achieve its leading constant (Corollary~\ref{thm:asymptotic}). The smooth and composite cases separate at second order, where only the smooth minimax rate depends on $L$ (last row of Table~\ref{tab:price}).\medskip

\noindent\textbf{Organization.\ }
Section~\ref{sec:classical} reviews the related work: the classical bounds for \algn{Prox-GD} and \algno, the exact deterministic rates, and the two last-iterate results of \citet{lz24} and \citet{gs26} to which ours are directly related. Section~\ref{sec:envelope} presents the Moreau-envelope approach: the facts about the envelope that we use, the deterministic step, the smooth case, where the result of \citet{gs26} is recovered, and the noisy proximal step, with the example showing why the smooth argument fails and the pointwise inequality that replaces it. Section~\ref{sec:main} contains the proof of Theorem~\ref{thm:exact}, with its attaining instance, and its consequences: the constant-stepsize corollary, the exact composite price, and the case $\ff=0$. 
Section~\ref{sec:design} is devoted to the design of optimal schedules and to their minimax optimality among all predetermined schedules, and Section~\ref{sec:asymptotic} to the asymptotic minimax law.
Appendix~\ref{app:certificate} gives the complete interpolation certificate behind Lemma~\ref{lem:local}, and Appendix~\ref{app:proofs} the proofs of Theorems~\ref{thm:convex}, \ref{thm:scalarconstruction}, and~\ref{thm:secondorder}.

\section{Related work}\label{sec:classical}

In this section, the stepsize is constant, $\gamma_t\equiv\gamma\in\big(0,\frac2L\big)$, and the variance bound is uniform, $\sigma_t\equiv\sigma \ge 0$, so that $S_t=(T-t)\gamma$. 
If there is no noise, that is, $e^t\equiv 0$, \algno reduces to \emph{Proximal Gradient Descent} (\algn{Prox-GD}), also known as the forward--backward algorithm 
\citep{cw05}.
We recall known bounds for \algn{Prox-GD} and \algno. 
All the bounds below are nondecreasing functions of $\|x^0-x^\star\|$, and we state them with $D$ in place of $\|x^0-x^\star\|$; a bound is said to be exact, or attained, when there is an instance with $\|x^0-x^\star\|=D$ on which it holds with equality. 

First, \algn{Prox-GD} has monotone objective gaps: by the descent lemma \citep[Lemma 10.4]{beck17},
$\Delta_{t+1}\leq \Delta_t +\left(\frac{L}{2}-\frac{1}{\gamma}\right)\big\|x^{t+1}-x^t\big\|^2\le\Delta_t$, 
for every $t\ge 0$, 
where the first inequality is tight on the quadratic $\ff=\frac L2\|\cdot\|^2$, $\RR=0$.
For \algn{Prox-GD} with $\gamma\le\frac1L$, the classical Lyapunov argument, which combines the squared distance to a solution and the objective gap, gives the rate \citep[Theorem 10.21]{beck17} 
\begin{align}
\Delta_t \le  \frac{D^2}{2\gamma t}, \quad\forall t\ge 1.\label{eqr1}
\end{align}
For comparison, when $\RR=0$ and $\gamma\le\frac1L$, the exact worst-case rate of Gradient Descent (\algn{GD}), first established by \citet{dt14}, is $\Delta_t\le\frac{D^2}{4\gamma t+2/L}$ \citep[Theorem 2]{tv23}; it is attained by a Huber function \citep[Theorem 3]{tv23},
so that, in this case, \eqref{eqr1} is suboptimal, asymptotically by a factor of 2.
For $\gamma\in\big(\frac1L,\frac2L\big)$, the exact worst-case rate of \algn{GD}, conjectured by \citet{dt14} and proved by \citet{kim24}, is
\begin{equation*}
\Delta_t\le\frac{LD^2}{2}\max\bigg\{\frac{1}{1+2\gamma Lt},(1-\gamma L)^{2t}\bigg\},\quad\forall t\ge1.
\end{equation*}
The second term, attained by the quadratic $\ff=\frac L2\|\cdot\|^2$, is the larger one if and only if $\gamma L>\beta_t$, where $\beta_t$ denotes the unique root in $(1,2)$ of the equation $(1-\beta)^{-2t}-1=2t\beta$ \citep[Definition 2.11]{rgp25}; the sequence $(\beta_t)_{t\ge1}$ is increasing, with $\beta_1=\frac32$ and $\beta_t\to2$. Hence, the Huber term $\frac{D^2}{4\gamma t+2/L}$ is the exact worst case of \algn{GD} for every $t\ge1$ if and only if $\gamma L\le\beta_1=\frac32$. This is the origin of the stepsize limit $\frac{3}{2L}$ 
that appears in our analysis through the admissibility condition \eqref{eq:admissible}, which reads $L\gamma_{T-1}\le\beta_1$ at the last iteration, where $S_{T-1}=\gamma_{T-1}$.

In the composite case, the exact worst-case rate of \algn{Prox-GD} for $\gamma\le\frac1L$ is $\Delta_t\le\frac{D^2}{4\gamma t}$ \citep[Theorems 7 and 8]{tv23};
it is already attained in the case $\ff=0$, which is $L$-smooth for every $L>0$, with $\RR=\lambda|\cdot|$ in dimension one \citep[Theorem 4.1]{thg17}. Thus, in the composite case, \eqref{eqr1} is suboptimal by a factor of exactly 2.
For $\gamma\in(\frac1L,\frac2L)$, the two examples just mentioned show that the worst case of \algn{Prox-GD} is at least $\max\big\{\frac{D^2}{4\gamma t},\frac{LD^2}{2}(1-\gamma L)^{2t}\big\}$, but its exact value was not known beyond $\frac1L$; Corollary~\ref{cor:constant} below shows that it is $\frac{D^2}{4\gamma t}$, for $\gamma$ up to $\frac{3}{2L}$.
Let us also mention that \algn{Prox-GD} can be accelerated by nonconstant schedules with large stepsizes \citep{ba25}, 
so that an optimal schedule need not have small stepsizes; this is why the design problem of Section~\ref{sec:design} and the lower bounds of Section~\ref{sec:asymptotic} are posed over all positive schedules, with no upper bound on the stepsizes.
For \algn{GD}, \citet{gri24} proved faster worst-case guarantees with periodic stepsizes much larger than $\frac2L$, and left open their effect in the proximal and stochastic settings; in our setting, the minimax optimal schedule has no stepsize larger than $\frac{3}{2L}$, unless the noise level is low (Theorem~\ref{thm:global} and Corollary~\ref{cor:globalT}).

When $\sigma>0$ and $\gamma\le\frac1L$, the same kind of Lyapunov argument
yields a 
bound on the averaged iterate $\bar{x}^T\eqdef\frac1T\sum_{t=1}^Tx^t$: by convexity of $\Psi$,
\begin{equation}
\E\big[\Psi(\bar{x}^T)\big]-\Psi^\star\le\frac1T\sum_{t=1}^T\E[\Delta_t]\le\frac{D^2}{2\gamma T}+\gamma\sigma^2,\quad\forall T\ge1.\label{eqavg}
\end{equation}
However, this argument says nothing about the last iterate, since $(\Delta_t)_{t\ge0}$ is no longer monotone when $\sigma>0$.
In this setting, last-iterate bounds for \algno have been obtained only recently: \citet[Lemma 4.2]{lz24} proved that, for $\gamma\le\frac{1}{2L}$,
\begin{equation}
\E[\Delta_T]\le\frac{D^2}{2\gamma T}+2\gamma\sigma^2\sum_{t=1}^{T}\frac1t,\quad\forall T\ge1.\label{eq:lz}
\end{equation}
This is the specialization to our setting of a bound valid for composite mirror descent, in a non-Euclidean geometry; 
they complement it with bounds of the optimal order $O\big(\frac{LD^2}{T}+\frac{\sigma D}{\sqrt T}\big)$ for linearly decaying stepsizes \citep[Theorems 3.4 and 3.5]{lz24}. The initial-error term in \eqref{eq:lz} is that of \eqref{eqr1} and \eqref{eqavg}, and the noise term grows like $2\gamma\sigma^2\log T$.

In the smooth case $\RR=0$, \citet[Theorem 1.1]{gs26} obtained the exact worst case of \algn{SGD}: for $\gamma\le\frac{3}{2L}$,
\begin{equation}
\E[\Delta_T]\le\frac{D^2}{4\gamma T+2/L}+\frac{L\gamma^2\sigma^2}{2}\sum_{k=0}^{T-1}\frac{1}{1+2kL\gamma},\quad\forall T\ge1,\label{eq:gs}
\end{equation}
and the bound is attained 
by a separable instance in dimension $T+1$ made of Huber functions with Rademacher noise. More generally, their result gives the exact worst case of \algn{SGD} for every stepsize schedule in $\big[0,\frac{3}{2L}\big]^T$: the upper bound holds in any dimension and the lower bound holds for every nonnegative schedule in dimension at least $T+1$. 
The first term in \eqref{eq:gs} is the exact deterministic rate recalled above, and the range $\gamma L\le\frac32=\beta_1$ is the one for which that rate is exact for every $T\ge1$. The second term is the sum, over the injections, of the variance $\gamma^2\sigma^2$ of the perturbation weighted by half the curvature $\frac{L}{1+2kL\gamma}$ felt by an injection followed by $k$ more iterations;
it grows like $\frac{\gamma\sigma^2}{4}\log T$, so that the logarithm is unavoidable with a constant stepsize.
Compared with \eqref{eq:lz}, the constant of the initial-error term is improved by a factor of 2 and the constant in front of the logarithm by a factor of 8. Minimizing the exact formula over the schedule, \citet{gs26} further obtained the minimax optimal schedule of \algn{SGD},
provided that $\sigma^2\ge\frac{2}{3\sqrt3}\frac{L^2D^2}{\sqrt T}$: it is decreasing, close to linearly decaying, and achieves the rate $\frac{\sigma D}{\sqrt T}-\frac{\sigma^{3/2}\sqrt D}{\sqrt L\,T^{3/4}}+o(T^{-3/4})$,
without a logarithm.
Their proof is an induction on an inequality involving an auxiliary sequence of points, convex combinations of $x^\star$ and of the points $x^s-\frac1L\nabla\ff(x^s)$, $s<t$, and the interpolation inequalities between all the iterates. 

Last-iterate guarantees have also been studied for subgradient and stochastic subgradient methods on nonsmooth objectives, under a bound on the subgradients \citep{sz13,jnn19,dcmm24}, including the exact worst-case value and the optimal schedule in the deterministic case \citep{zg25}. This is another problem class, with another oracle: the nonsmooth part is accessed through its subgradients and their magnitude governs the rate, whereas in our setting it is handled exactly by its prox and only the gradient of the smooth part is noisy. These results are not comparable with ours.

In the composite case, no exact last-iterate result was available, and \eqref{eq:lz} and \eqref{eq:gs} leave open the price of the prox, both on the initial-error term, whose exact value was known for $\gamma\le\frac1L$ only, and on the noise term, where the weight $\frac{4}{k+1}$ of an injection followed by $k$ more iterations in \eqref{eq:lz} and the smooth weight $\frac{L\gamma}{1+2kL\gamma}$ in \eqref{eq:gs}, both relative to $\frac{\gamma\sigma^2}{2}$,
differ by a factor of 8 asymptotically and by much more for the last injections when $L\gamma$ is small.
We show in the next two sections that the envelope potential extends to the composite case, provided that the whole noisy proximal step is controlled at once, and that the resulting bound is exact, for every stepsize $\gamma\le\frac{3}{2L}$.

\section{The Moreau-envelope approach}\label{sec:envelope}

Our analysis rests on a potential that is the \emph{Moreau envelope} of $\Psi$ at the current iterate $x^t$, with a smoothing parameter that decreases along the run and vanishes at the horizon: $\tau_t\eqdef2S_t$, where $S_t$ is the remaining stepsize sum \eqref{eq:tails}, so that $\tau_t=2(T-t)\gamma$ for a constant stepsize. Each step removes $2\gamma_t$ from the parameter and $\tau_T=0$, where the envelope is $\Psi$ itself, so that a telescoping estimate on this potential ends at the objective gap of the last iterate, and not at an average of earlier gaps as in \eqref{eqavg}. The reason why this potential handles the noise correctly is its curvature. The noise injected at iteration $t$ is felt by a potential through its curvature, which is $L$ for $\Psi$ itself when $\RR=0$, and only $L/(1+L\tau)$ for the envelope of $\ff$ with parameter $\tau$ (Lemma~\ref{lemenvsmooth} below): the envelope is nearly flat while the horizon is far, and it stiffens progressively, up to $\Psi$ at $t=T$. 
As a result, with a constant stepsize, a noise injection followed by $k$ more iterations is damped by a factor of order $1/(1+2kL\gamma)$, which is exactly the worst-case decay of the objective gap of \algn{GD} over $k$ iterations, recalled in Section~\ref{sec:classical}, and the sum of these factors over the injections grows only logarithmically with $T$.
The potential depends on the horizon, but the algorithm does not: with a constant stepsize, the bounds obtained hold for every $T\ge1$ along a single run; with a horizon-dependent schedule, the horizon is the one the schedule was designed for.

The idea of using the Moreau envelope is not new: \citet[Section 2.1]{gs26} observed that the inductive hypothesis of their analysis of \algn{SGD} yields a bound on the expected value of the envelope of $\ff$ at the current iterate, with the same parameter $\tau_t$. Their induction, however, is not on this quantity: it is on a stronger hypothesis, which involves an auxiliary sequence built from all the past iterates, and from which the bound on the envelope is derived.
What is new is to make the envelope the potential itself: we show that a bound on it propagates from one iteration to the next by itself, without the history of the iterates. This is what allows us to handle the composite case $\RR\neq0$, for which we also provide a matching lower bound.
We first collect the facts about the envelope that we use (Section~\ref{sec:envelopefacts}); we then show that a deterministic step does not increase the envelope, once its parameter is decreased by twice the stepsize (Section~\ref{sec:detstep}); the smooth case follows at once (Section~\ref{sec:smoothcase});
finally, we explain why this argument breaks down when $\RR\neq0$ and establish
the pointwise inequality that controls the whole noisy proximal step (Section~\ref{sec:noisystep}).

\subsection{The Moreau envelope}\label{sec:envelopefacts}
For every $\tau>0$, define the Moreau envelope and proximal point of the full objective: for every $x\in\cH$,
\begin{equation}
 \Psi_\tau(x)\eqdef\min_{y\in\cH}\left\{\Psi(y)+\tfrac1{2\tau}\norm{y-x}^2\right\},
 \qquad p\eqdef\prox_{\tau\Psi}(x),\qquad \Psi_0\eqdef\Psi.\label{eq:envelope}
\end{equation}
For $\tau>0$, the minimum in \eqref{eq:envelope} is attained at the unique point $p=\prox_{\tau\Psi}(x)$, which is characterized by 
$x-p=\tau\big(\nabla \ff(p)+s_p\big)$ for a unique $s_p\in\partial \RR(p)$ \citep[Proposition 16.44]{bau17}. We define $G_p\eqdef\nabla \ff(p)+s_p$. 
Evaluating \eqref{eq:envelope} at $y=p$ and comparing with $y=x^\star$, we obtain
\begin{align}
& G_p=\tfrac1{\tau}(x-p)\in\partial\Psi(p),\qquad
 \Psi_\tau(x)\overset{\eqref{eq:envelope}}{=}\Psi(p)+\tfrac\tau2\norm{G_p}^2,\label{eq:proxidentity}\\
 &0\leq\Psi_\tau(x)-\Psi^\star\overset{\eqref{eq:envelope}}{\leq} \tfrac1{2\tau}
                    \norm{x-x^\star}^2.\label{eq:initial}
\end{align}
Since $\prox_{\tau\Psi}(x^\star)=x^\star$, we have $\Psi_\tau(x^\star)=\Psi^\star$, so that $\Psi_\tau(x)-\Psi^\star$ in \eqref{eq:initial} is the objective gap of the envelope at $x$. The envelopes are a tool of the analysis only: the algorithm never evaluates $\prox_{\tau\Psi}$. When $\RR=0$, we write $\ff_\tau\eqdef\Psi_\tau$, and we have more:

\begin{lemma}\label{lemenvsmooth}
Suppose that $\RR=0$ and let $\tau\ge0$.
Then $\ff_\tau$ is convex and differentiable, and $\nabla \ff_\tau$ is $\frac{L}{1+L\tau}$-Lipschitz continuous.
Consequently, for every $x,e\in\cH$,
\begin{equation}
\ff_\tau(x+e)\le \ff_\tau(x)+\langle\nabla \ff_\tau(x),e\rangle+\frac{L}{2(1+L\tau)}\|e\|^2,\qquad
\ff_\tau(x)-\ff(x^\star)\le\frac{L}{2(1+L\tau)}\|x-x^\star\|^2.\label{eq:envsmooth}
\end{equation}
\end{lemma}
\begin{proof}
The case $\tau=0$ is the classical descent lemma, so let $\tau>0$. 
Convexity and differentiability of $\ff_\tau$, with $\nabla \ff_\tau(x)=\frac{x-p}{\tau}=\nabla \ff(p)$, where $p=\prox_{\tau \ff}(x)$, are classical \citep[Propositions 12.15 and 12.30]{bau17}.
Let $x,x'\in\cH$, and let $p,p'$ be their proximal points.
Then $x-x'=(p-p')+\tau\big(\nabla \ff(p)-\nabla \ff(p')\big)$ and, by the $\frac1L$-cocoercivity of $\nabla \ff$ \citep[Theorem 2.1.5]{nes18}, $\langle\nabla \ff_\tau(x)-\nabla \ff_\tau(x'),x-x'\rangle=\langle\nabla \ff(p)-\nabla \ff(p'),p-p'\rangle+\tau\|\nabla \ff(p)-\nabla \ff(p')\|^2\ge\big(\tau+\tfrac1L\big)\|\nabla \ff_\tau(x)-\nabla \ff_\tau(x')\|^2$.
Hence $\nabla \ff_\tau$ is $\big(\tau+\frac1L\big)$-cocoercive and, in particular, $\frac{L}{1+L\tau}$-Lipschitz continuous,
and the first inequality in \eqref{eq:envsmooth} is the descent lemma for $\ff_\tau$. The second one is the first one at $x^\star$, since $\nabla \ff_\tau(x^\star)=\nabla \ff(x^\star)=0$.
\end{proof}
The improved curvature $L/(1+L\tau)$ is specific to the smooth case. In general, for $\tau>0$, the envelope $\Psi_\tau$ is differentiable with $\nabla\Psi_\tau(x)=G_p=\frac{x-p}{\tau}$, which is $\frac1\tau$-Lipschitz continuous \citep[Proposition 12.30]{bau17}, and this is the best one can say: for $\ff=0$, which is $L$-smooth for every $L>0$, and $\RR=\lambda|\cdot|$ in dimension one, with $\lambda>0$, the envelope is $\Psi_\tau(z)=\frac{z^2}{2\tau}$ for $|z|\le\lambda\tau$, with curvature $\frac1\tau$ near zero, whereas $L/(1+L\tau)$ can be made arbitrarily small. The composite case thus loses the curvature of the envelope, both in the initial bound \eqref{eq:initial},
which does not see $L$, unlike the second inequality in \eqref{eq:envsmooth},
and, more seriously, in the treatment of the noise, as we will see in Section~\ref{sec:noisystep}.

\subsection{The deterministic step}\label{sec:detstep}
The key property is that a step of \algn{Prox-GD} does not increase the Moreau envelope, provided that its parameter is decreased by $2\gamma$. We use the interpolation inequality for convex $L$-smooth functions \citep[Theorem 2.1.5]{nes18}: for every $y,z\in\cH$,
\begin{equation}
\ff(y)\ge \ff(z)+\langle\nabla \ff(z),y-z\rangle+\tfrac{1}{2L}\|\nabla \ff(y)-\nabla \ff(z)\|^2,\label{eq:interp}
\end{equation}
and, for $\RR$, the subgradient inequality: for every $y,z\in\cH$ and $s_z\in\partial \RR(z)$,
\begin{equation}
\RR(y)\ge \RR(z)+\langle s_z,y-z\rangle.\label{eq:subgradient}
\end{equation}

\begin{lemma}\label{lemdet}
Let $\gamma>0$ and $\tau\ge2\gamma$. Set $\theta\eqdef\gamma/\tau\in\big(0,\frac12\big]$ and suppose that $\theta+L\gamma\le2$. Then, for every $x\in\cH$,
\begin{equation}
\Psi_{\tau-2\gamma}(q)\le\Psi_\tau(x),\qquad\text{where}\quad q\eqdef\prox_{\gamma \RR}\big(x-\gamma\nabla \ff(x)\big).\label{eq:det}
\end{equation}
\end{lemma}
\begin{proof}
Set $\tau'\eqdef\tau-2\gamma\ge0$, $p\eqdef\prox_{\tau\Psi}(x)$, and $u\eqdef\prox_{\tau'\Psi}(q)$
($u\eqdef q$ if $\tau'=0$).
By the characterization of the proximal points, there exist
$s_p\in\partial \RR(p)$, $s_q\in\partial \RR(q)$, $s_u\in\partial \RR(u)$ (with $s_u\eqdef s_q$ if $\tau'=0$)
such that, with the notation $g_x\eqdef\nabla \ff(x)$,
$g_p\eqdef\nabla \ff(p)$, $g_u\eqdef\nabla \ff(u)$, $G_p\eqdef g_p+s_p$, $G_x\eqdef g_x+s_q$, $U\eqdef g_u+s_u$,
\begin{equation}
x-p=\tau G_p,\qquad x-q=\gamma G_x,\qquad q-u=\tau'U.\label{eq:geom}
\end{equation}
By \eqref{eq:proxidentity},
\begin{equation}
\Psi_\tau(x)-\Psi_{\tau'}(q)\overset{\eqref{eq:proxidentity}}{=}\ff(p)-\ff(u)+\RR(p)-\RR(u)+\tfrac{\tau}{2}\|G_p\|^2-\tfrac{\tau'}{2}\|U\|^2.\label{eq:diff}
\end{equation}
We lower bound $\ff(p)-\ff(u)$ by going from $u$ to $p$ partly directly and partly through $x$, writing $\ff(p)-\ff(u)=(1-\theta)\big[\ff(p)-\ff(u)\big]+\theta\big[\ff(x)-\ff(u)\big]+\theta\big[\ff(p)-\ff(x)\big]$ and applying \eqref{eq:interp} to each bracket:
\begin{align}
\ff(p)-\ff(u)\overset{\eqref{eq:interp}}{\ge}{}&(1-\theta)\Big[\langle g_u,p-u\rangle+\tfrac{1}{2L}\|g_p-g_u\|^2\Big]+\theta\Big[\langle g_u,x-u\rangle+\tfrac{1}{2L}\|g_x-g_u\|^2\Big]\notag\\
&+\theta\Big[\langle g_x,p-x\rangle+\tfrac{1}{2L}\|g_p-g_x\|^2\Big].\label{eq:weightedf}
\end{align}
Similarly, for $\RR$, we go from $u$ to $p$ partly directly and partly through $q$, the point at which a subgradient of $\RR$ is available, with the weights $1-h,h,h$, where $h\eqdef\frac{\gamma}{\tau-\gamma}\in(0,1]$:
\begin{equation}
\RR(p)-\RR(u)\overset{\eqref{eq:subgradient}}{\ge}(1-h)\langle s_u,p-u\rangle+h\langle s_u,q-u\rangle+h\langle s_q,p-q\rangle.\label{eq:weightedR}
\end{equation}
All the inner products are evaluated using \eqref{eq:geom}: $p-u=\gamma G_x+\tau' U-\tau G_p$, $x-u=\gamma G_x+\tau' U$, $p-x=-\tau G_p$, $q-u=\tau' U$, $p-q=\gamma G_x-\tau G_p$. Plugging the two lower bounds into \eqref{eq:diff} yields a quadratic form in the six vectors $g_x,g_p,g_u,s_p,s_q,s_u$, which, in terms of
\begin{equation*}
a\eqdef G_p-U,\qquad b\eqdef G_x-U,\qquad A_f\eqdef g_p-g_u,\qquad B_f\eqdef g_x-g_u\qquad(\text{so that } b-B_f=s_q-s_u),
\end{equation*}
reads
\begin{align}
\Psi_\tau(x)-\Psi_{\tau'}(q)\overset{\eqref{eq:diff},\eqref{eq:geom},\eqref{eq:weightedf},\eqref{eq:weightedR}}{\ge}{}&\frac{\tau}{2}\|a\|^2-\frac{\gamma\tau}{\tau-\gamma}\langle a,b\rangle+\frac{\gamma^2}{\tau-\gamma}\|b\|^2+\frac{\gamma^2}{\tau-\gamma}\langle a-b,B_f\rangle\notag\\
&+\frac{1}{2L}\|A_f\|^2-\frac{\gamma}{L\tau}\langle A_f,B_f\rangle+\frac{\gamma}{L\tau}\|B_f\|^2\notag\\
={}&\frac{\tau}{2}\Big\|a-h\big(b-\theta B_f\big)\Big\|^2
+\frac{1}{2L}\big\|A_f-\theta B_f\big\|^2\notag\\
&+\frac{\gamma^2(\tau-2\gamma)}{2(\tau-\gamma)^2}\|s_q-s_u\|^2
+\frac{\theta(2-\theta-L\gamma)}{2L}\|B_f\|^2,\label{eq:sos}
\end{align}
as one checks by expanding the squares. Every coefficient in \eqref{eq:sos} is nonnegative under the assumptions $\tau\ge2\gamma$ and $\theta+L\gamma\le2$, which proves \eqref{eq:det}.
\end{proof}

\noindent\textbf{How the weights are found.\ }
The weights $\theta$ and $h$ in the proof are not guessed: they are forced by cancellation. Compare $\ff(p)$ with $\ff(u)$ either directly or through $x$, and $\RR(p)$ with $\RR(u)$ either directly or through $q$, the point at which a subgradient of $\RR$ is available, and give the two indirect paths weights $\vartheta,\varkappa\in[0,1]$, to be determined: the three interpolation inequalities \eqref{eq:interp} receive the weights $(1-\vartheta,\vartheta,\vartheta)$ and the three subgradient inequalities the weights $(1-\varkappa,\varkappa,\varkappa)$. The function values telescope for any weights, and the difference \eqref{eq:diff}, minus the six weighted interpolation residuals (which are written explicitly in Appendix~\ref{app:certificate}), is a quadratic form in the six vectors $g_x,g_p,g_u,s_p,s_q,s_u$. The certificate we look for is a sum of squares of the four differences $a,b,A_f,B_f$, in which the composite gradient $U$ at $u$ does not appear, so that the terms involving $U$ must cancel. After substitution of the displacement relations \eqref{eq:geom}, these terms are
\begin{equation}
 [\gamma-\varkappa(\tau-\gamma)]\ip Ub
 +[\varkappa(\tau-\gamma)-\vartheta\tau]\ip U{B_f},
 \label{eq:discoverycancel}
\end{equation}
the coefficients of $\norm U^2$ and $\ip Ua$ canceling already, for any weights. Removing the two cross terms in \eqref{eq:discoverycancel} forces $
 \varkappa=\tfrac\gamma{\tau-\gamma}=h$, 
 $\vartheta=\tfrac\gamma\tau=\theta$.
The six weights are then nonnegative if and only if $\tau\geq2\gamma$, and the resulting quadratic form is the one in \eqref{eq:sos}; completing its squares exposes the other restriction $\theta+L\gamma\leq2$. When the noise is restored, in Lemma~\ref{lem:local} below, $q$ becomes the noisy output of the step, $G_x$ becomes $g_x+e+s_q$, and the only change in the algebra is the nonsmooth displacement, $s_q-s_u=b-B_f-e$, which brings in the term $\gamma\ip Ue$. This term is canceled by the term $-\gamma G_p=-\gamma(U+a)$ of the vector $B(x)$ of the lemma, given in \eqref{eq:Bc}, and the four squares of \eqref{eq:sos} are recentered by multiples of $e$; the exact identity \eqref{eq:certificate} of Appendix~\ref{app:certificate} gives the second term of $B(x)$, a multiple of the gradient difference $A_f-B_f=\nabla\ff(p)-\nabla\ff(x)$, and the noise coefficient $c$. 
Thus the weights come from cancellation, the admissibility condition from positivity, and the noise coefficient from completing the same squares. This is the whole certificate when $L\gamma\le1$. When $L\gamma>1$, the noise coefficient is improved by involving one more point, the noiseless output $\qc\eqdef\prox_{\gamma\RR}\big(x-\gamma\nabla\ff(x)\big)$ of the step, through three further subgradient inequalities, between $q$ and $\qc$ in both directions and between $\qc$ and $u$, with weights proportional to $L\gamma-1$; this adds a third term to $B(x)$ and a fifth square (Appendix~\ref{app:long}). Proposition~\ref{prop:localsharp} shows that the resulting coefficient is unavoidable, in both cases.

\begin{remark}\label{remdet}
If $\RR=0$, then $s_p=s_q=s_u=0$, $a=A_f$, $b=B_f$, and since $h(1-\theta)=\theta$, the first two squares in \eqref{eq:sos} add up to $\frac{1+L\tau}{2L}\|a-\theta b\|^2$, while the third one vanishes. The third square is the specific ingredient of the composite case: the subgradients of $\RR$ at $q$ and at $u$ may differ, and its coefficient vanishes when $\tau=2\gamma$. The two restrictions of the lemma are sharp. 
First, $2\gamma$ is the largest uniform decrement of the parameter, in the smooth case as well: if $\ff=0$ and $\RR=\lambda|\cdot|$ in dimension one, with $\lambda>0$, then $\Psi_\tau(x)=\lambda x-\frac{\lambda^2\tau}{2}$ for $x\ge\lambda\tau$ and $q=x-\gamma\lambda$, so that, for $x$ large enough and every $\nu\in[0,\tau]$, $\Psi_{\tau-\nu}(q)-\Psi_\tau(x)=\lambda^2\big(\frac\nu2-\gamma\big)$: the envelope does not increase if and only if $\nu\le2\gamma$, and the same computation holds on the linear facet of a Huber function with $\RR=0$.
Smaller decrements are allowed, since $\Psi_\tau(x)$ is nonincreasing in $\tau$. Second, the condition $\theta+L\gamma\le 2$ is sharp: if $\RR=0$ and $\ff=\frac L2\|\cdot\|^2$, then $\ff_\tau(x)=\frac{L}{2(1+L\tau)}\|x\|^2$ and $q=(1-\gamma L)x$, so that \eqref{eq:det} holds if and only if $(1-\gamma L)^2(1+L\tau)\le1+L\tau-2L\gamma$, which is equivalent to $\theta+L\gamma\le2$.
\end{remark}

\subsection{The smooth case}\label{sec:smoothcase}

When $\RR=0$, the noise is handled by the smoothness of the envelope, Lemma~\ref{lemenvsmooth}, and we recover the bound of \citet{gs26}.

\begin{proposition}[The smooth case]\label{prop:smooth}
If $\RR=0$ and the schedule $\boldsymbol{\gamma}$ is admissible, that is, satisfies \eqref{eq:admissible}, then
\begin{equation}
 \E[\ff(x^T)-\ff(x^\star)]\leq
 F_{L,T}^{\mathrm{smooth}}(\boldsymbol{\gamma};D,\boldsymbol{\sigma})\eqdef 
 \frac{D^2}{4S_0+2/L}
 +\frac12\sum_{t=0}^{T-1}\frac{\sigma_t^2\gamma_t^2}{2S_{t+1}+1/L}.\label{eq:smooth}
\end{equation}
\end{proposition}
\begin{proof}
Set $\xc^{t+1}\eqdef x^t-\gamma_t\nabla \ff(x^t)$, which is $\cF^t$-measurable, so that $x^{t+1}=\xc^{t+1}-\gamma_te^t$
(all the expectations below are finite; see the integrability paragraph of Section~\ref{sec:main}).
Lemma~\ref{lemdet} with $\tau=2S_t$ and $\gamma=\gamma_t$, whose condition $\theta+L\gamma\le2$ is \eqref{eq:admissible}, gives $\ff_{2S_{t+1}}(\xc^{t+1})\le \ff_{2S_t}(x^t)$. The first inequality in \eqref{eq:envsmooth}, applied to $\ff_{2S_{t+1}}$ at $\xc^{t+1}$ with $e=-\gamma_te^t$, and the conditional unbiasedness of $e^t$ then give, after conditional and then total expectation,
\begin{equation*}
 \E[\ff_{2S_{t+1}}(x^{t+1})]
 \overset{\eqref{eq:envsmooth},\eqref{eq:det},\eqref{eq:noise}}{\leq}
 \E[\ff_{2S_t}(x^t)]+\frac{L\gamma_t^2\sigma_t^2}{2(1+2LS_{t+1})}.
\end{equation*}
Summing over $t=0,\ldots,T-1$, with $S_T=0$, and applying the second inequality in \eqref{eq:envsmooth} at $x^0$ with $\tau=2S_0$, we obtain \eqref{eq:smooth}.
\end{proof}
With a constant stepsize $\gamma_t\equiv\gamma$ and a uniform variance bound $\sigma_t\equiv\sigma$, we have $S_t=(T-t)\gamma$, and \eqref{eq:smooth} is the bound
\eqref{eq:gs} of \citet[Theorem~1.1]{gs26} recalled in Section~\ref{sec:classical}; the condition \eqref{eq:admissible} is weaker than $\gamma_t\le\frac{3}{2L}$ early in the run, when the remaining sum $S_t$ is large, and the two conditions coincide at the last iteration.
Compared with their induction, the proof uses only three interpolation inequalities, between $x$ and the two proximal points $p$ and $u$, and no auxiliary sequence;
this is what allows the extension to $\RR\neq0$ in the next subsection, where the gradient steps of their auxiliary sequence have no counterpart. When $\sigma_t\equiv0$, \eqref{eq:smooth} is the bound $\frac{D^2}{4S_0+2/L}$,
which, for a constant stepsize, is the exact rate of \algn{GD} recalled in Section~\ref{sec:classical}.
With a constant stepsize, the bound \eqref{eq:gs} also has a neat superposition structure: writing $\mathrm{GD}_k(D)\eqdef\frac{LD^2}{2(1+2kL\gamma)}$ for the exact worst case of $\Delta_k$ after $k$ iterations of \algn{GD} started at distance $D$ from $x^\star$, we have
\begin{equation*}
\E[\Delta_T]\overset{\eqref{eq:gs}}{\le}\mathrm{GD}_T(D)+{\textstyle\sum_{t=0}^{T-1}}\mathrm{GD}_{T-1-t}(\gamma\sigma):
\end{equation*}
every noise injection acts as a fresh initial distance $\gamma\sigma$, which the remaining iterations damp by the deterministic worst-case factor. In particular, the last injection is not damped at all, which gives the term $\frac{L\gamma^2\sigma^2}{2}$, and the sum of the damping factors is the logarithm.

\subsection{Why the smooth argument fails when \texorpdfstring{$\RR\neq0$}{R is nonzero}, and the noisy-step inequality}\label{sec:noisystep}
When $\RR\neq0$, the noise is \emph{inside} the prox, and the mapping $z\mapsto\Psi_{\tau'}\big(\prox_{\gamma\RR}(z)\big)$ is not smooth:
for $\RR=\lambda|\cdot|$ in dimension one, $\prox_{\gamma\RR}$ is the soft-thresholding operator $z\mapsto\operatorname{sign}(z)\max\{|z|-\gamma\lambda,0\}$, and at a point $z$ on the boundary of its dead zone $[-\gamma\lambda,\gamma\lambda]$, the two symmetric perturbations $z\pm\gamma\varepsilon$ are mapped to points at distances $\gamma\varepsilon$ and $0$ from $\prox_{\gamma\RR}(z)$. The displacement of the output does not average out, $\E\big[\prox_{\gamma\RR}(z-\gamma e)\big]\neq\prox_{\gamma\RR}(z)$ for a centered $e$, and if $\nabla\Psi_{\tau'}$ does not vanish at $\prox_{\gamma\RR}(z)$, the increase of the envelope is of first order in $\gamma\|e\|$, not of second order: there is no counterpart of \eqref{eq:envsmooth}. Here is a concrete instance of this failure, with $\tau=2$ and $\gamma=\frac12$, so that $\tau'=\tau-2\gamma=1$.
Take
\begin{equation}
 \ff(y)=\tfrac12(y-2)^2,\quad \RR(y)=|y|,\quad x=-3,\quad
 \gamma=\tfrac12,\quad e\in\{-\varepsilon,\varepsilon\},\quad 0<\varepsilon<2,
 \label{eq:biasexample}
\end{equation}
with equiprobable signs. The noiseless proximal output is $0$, whereas the noisy outputs are $0$ and $-\varepsilon/2$. Direct minimization gives
\begin{equation}
 \Psi_1(z)\overset{\eqref{eq:envelope},\eqref{eq:biasexample}}{=}
 \tfrac74-\tfrac12z+\tfrac14z^2,\qquad -1<z\leq0.
 \label{eq:biasenvelope}
\end{equation}
Consequently, for the noisy output $q_e$,
$\E[\Psi_1(q_e)]-\Psi_1(0)=\tfrac{\varepsilon}{8}+\tfrac{\varepsilon^2}{32}$.
The expected increase has a term linear in $\varepsilon$, which rules out bounding it by a constant times the variance $\varepsilon^2$ of the perturbation, as \eqref{eq:envsmooth} does when $\RR=0$.
What saves the argument is that this first-order loss occurs only when the inequality \eqref{eq:det} is strict,
and that the slack pays for it: in the example, $\Psi_2(-3)=\frac{17}{4}$, whereas $\Psi_1(0)=\frac74$. The right formulation is therefore a pointwise inequality for the whole step, from $x$ to the noisy output, in which the first-order term in the noise is linear, with a coefficient that does not depend on the noise, so that it disappears in conditional expectation, and in which the coefficient of the quadratic term is the price of the injection.

\begin{lemma}[Predictable noisy-step inequality]\label{lem:local}
Let $\gamma>0$, $\tau\geq2\gamma$, $\theta\eqdef\gamma/\tau$, and suppose that $\theta+L\gamma\leq2$. Set $\tau'\eqdef\tau-2\gamma$. There is a Lipschitz continuous mapping $B:\cH\to\cH$, which depends on $\ff$, $\RR$, $L$, $\gamma$, and $\tau$ only, such that, for every $x,e\in\cH$,
\begin{equation}
 \Psi_{\tau'}(q)\leq\Psi_\tau(x)+\ip{B(x)}e+c\norm e^2,\qquad\text{where}\quad
 q\eqdef\prox_{\gamma \RR}\big(x-\gamma(\nabla \ff(x)+e)\big)\label{eq:local}
\end{equation}
and
\begin{equation}
 c\eqdef\frac{\gamma^2}{2d},\qquad
 d\eqdef\tau'+\omega(\gamma)\overset{\eqref{eq:omega}}{=}\begin{cases}\tau-L\gamma^2,&\text{if }L\gamma\le1,\\[0.5mm] \tau'+\frac1L,&\text{if }L\gamma\ge1.\end{cases}\label{eq:c}
\end{equation}
The denominator $d$ is positive. The expression of $B(x)$ is given in \eqref{eq:Bc}, in Appendix~\ref{app:certificate}.
\end{lemma}
\begin{proof}
Appendix~\ref{app:certificate} gives self-contained identities expressing the difference between the right-hand and left-hand sides of \eqref{eq:local}, with $B(x)$ defined in \eqref{eq:Bc}, as a sum of nonnegative interpolation residuals and squared norms: six residuals and four squares when $L\gamma\le1$, in \eqref{eq:certificate}, and nine residuals and five squares when $L\gamma>1$, in \eqref{eq:certificatelong}.
All coefficients are nonnegative under the stated hypotheses, including at $\tau'=0$. This proves the pointwise inequality without any probabilistic assumptions.
\end{proof}
The choice $e=0$ gives Lemma~\ref{lemdet} again, and the two restrictions on $\gamma$ and $\tau$ are the same, so that the sharpness discussion of Remark~\ref{remdet} applies. The vector $B(x)$ does not depend on $e$, which is what makes the linear term harmless: when the lemma is applied at $x=x^t$ with $e=e^t$, $B(x^t)$ is $\cF^t$-measurable and the conditional expectation of $\ip{B(x^t)}{e^t}$ vanishes. 
The noise coefficient $c$ has two regimes. We say that the step is \emph{short} if $L\gamma\le1$ and \emph{long} otherwise: for a short step, $c$ is $\frac{\gamma^2}{2(\tau-L\gamma^2)}$; for a long step, it is $\frac{L\gamma^2}{2(1+L\tau')}$, the coefficient of the smooth case in \eqref{eq:envsmooth}.
The following proposition shows that $c$ is the smallest possible uniform coefficient, in both regimes. Its attaining instance, like those of Section~\ref{sec:main},
is built from the scalar Huber function with threshold $\zeta>0$,
\begin{equation}
 \hub_\zeta(z)\eqdef\begin{cases}Lz^2/2,&|z|\leq\zeta,\\L\zeta|z|-L\zeta^2/2,&|z|\geq\zeta,\end{cases}\label{eq:huber}
\end{equation}
which is convex, with the $L$-Lipschitz derivative $\hub_\zeta'(z)=L\max\{-\zeta,\min\{z,\zeta\}\}$: quadratic on $[-\zeta,\zeta]$, it is linear with slope $\pm L\zeta$ outside.

\begin{proposition}[Optimality of the noise coefficient]\label{prop:localsharp}
Let $\gamma$ and $\tau$ satisfy the hypotheses of Lemma~\ref{lem:local}, and let $\sigma>0$.
There exist a convex function $\ff$ with $L$-Lipschitz derivative and a convex function $\RR$, both defined on $\R$, and a point $x\in\R$ such that, with the notation of Lemma~\ref{lem:local},
\begin{equation}
 \Psi_{\tau'}(q)-\Psi_\tau(x)=c\,\sigma^2\qquad\text{for $e=\sigma$ and for $e=-\sigma$.}\label{eq:localequality}
\end{equation}
Consequently, the coefficient $c$ in \eqref{eq:c} cannot be decreased in an inequality of the form \eqref{eq:local} valid for all functions $\ff$ and $\RR$ and all $x$ and $e$, whatever the vector $B(x)$, as long as it does not depend on $e$.
\end{proposition}
\begin{proof}
Set $\gc\eqdef\min\big\{\gamma,\frac1L\big\}$, so that $d=\tau'+\gc(2-L\gc)$, and
\begin{equation}
 w\eqdef\tfrac{\gamma\sigma}{d},\qquad
 \zeta\eqdef\gc w,\qquad \lambda\eqdef(1-L\gc)w\geq0,\qquad
 \ff\eqdef\hub_\zeta,\quad \RR\eqdef\lambda|\cdot|,\quad x\eqdef0.\label{eq:localextremal}
\end{equation}
Let $e\in\{-\sigma,\sigma\}$. Soft thresholding gives $q=-\operatorname{sign}(e)\gamma(\sigma-\lambda)=-\operatorname{sign}(e)(\zeta+\tau'w)$; the second equality holds because $\gamma\sigma=dw$ and $(\gamma-\gc)(1-L\gc)=0$.
At $u\eqdef-\operatorname{sign}(e)\zeta$, the derivative of $\Psi$ is $-\operatorname{sign}(e)w$, so that $q-u$ is $\tau'$ times that derivative. Thus $u=\prox_{\tau'\Psi}(q)$ when $\tau'>0$, and $u=q$ when $\tau'=0$. Consequently,
\begin{align*}
 \Psi_{\tau'}(q)
 &\overset{\eqref{eq:envelope},\eqref{eq:huber},\eqref{eq:localextremal}}{=}
 \tfrac L2\zeta^2+\lambda \zeta+\tfrac{\tau'}2w^2\overset{\eqref{eq:localextremal},\eqref{eq:c}}{=}c\sigma^2.
\end{align*}
Since $\Psi_\tau(0)=0$, this is \eqref{eq:localequality}. If an inequality of the form \eqref{eq:local} holds with a coefficient $c'$ and a vector $B(x)$ that does not depend on $e$, averaging it over the two signs of $e$ gives $c\sigma^2\le c'\sigma^2$.
\end{proof}

The form of the function $\omega$ in \eqref{eq:omega}, with its threshold $\frac1L$,
can be understood by considering a single noisy step from a solution.
Let $x=x^\star$ and $\tau=2\gamma$, with $L\gamma\le\frac32$: then $\tau'=0$, the vector $B(x^\star)$ of \eqref{eq:Bc} is zero, since $p=\qc=x^\star$ there, and \eqref{eq:local} reads $\Psi(q)-\Psi^\star\le\frac{\gamma^2}{2\omega(\gamma)}\norm e^2$. In dimension one, with $e=\pm\sigma$,
two mechanisms cause a loss at such a step. \emph{Curvature}: with $\ff=\frac L2|\cdot|^2$ and $\RR=0$, we have $q=-\gamma e$ and the loss is $\frac L2\gamma^2\sigma^2$.
\emph{A kink}: with $\ff=0$ and $\RR=\lambda|\cdot|$, $0\le\lambda\le\sigma$, soft thresholding gives $|q|=\gamma(\sigma-\lambda)$ and the loss is $\lambda\gamma(\sigma-\lambda)$, at most $\frac{\gamma\sigma^2}{4}$. The instance \eqref{eq:localextremal} combines them: for a short step, its slope $w=L\zeta+\lambda$ at the point reached is made of curvature and kink in the proportions $L\gamma$ and $1-L\gamma$, and its loss $\frac{\gamma\sigma^2}{2(2-L\gamma)}=\frac{\gamma^2\sigma^2}{2\omega(\gamma)}$ is larger than the loss of either mechanism alone. At $L\gamma=1$, the kink has weight zero. For a long step, it would need a negative weight, which is impossible, and curvature alone remains, with the loss $\frac L2\gamma^2\sigma^2$; this is the value $\omega(\gamma)=\frac1L$.

When stepsizes remain after the step, that is, $\tau'>0$, the comparison is the same:
since $\omega(\gamma)\le\frac1L$, with equality if and only if $L\gamma\ge1$,
\begin{equation*}
c\overset{\eqref{eq:c}}{=}\frac{\gamma^2}{2\big(\tau'+\omega(\gamma)\big)}\ \overset{\eqref{eq:omega}}{\ge}\ \frac{\gamma^2}{2\big(\tau'+1/L\big)}=\frac{L\gamma^2}{2(1+L\tau')},
\end{equation*}
which is the smooth constant of \eqref{eq:envsmooth}. In the case $\ff=0$, which is $L$-smooth for every $L>0$, or equivalently in the limit $L\to0$, $c=\frac{\gamma^2}{2\tau}$, whereas the smooth constant vanishes: a noise injection costs $\frac{\gamma^2\sigma^2}{2\tau}$ even without curvature, because the prox does not average out the noise. Proposition~\ref{prop:localsharp} shows that the price $c\sigma^2$ is paid by an actual instance at a single step; Theorem~\ref{thm:exact}, proved in the next section, shows that it accumulates along the whole run into the exact worst-case value of $\E[\Delta_T]$.

\section{Proof of the main theorem and its consequences}\label{sec:main}
We now prove Theorem~\ref{thm:exact}, stated in Section~\ref{sec:intro}, and derive its consequences. The proof has three parts: the integrability of the quantities involved, the upper bound, which is the telescoping of the noisy-step inequality of Lemma~\ref{lem:local} along the potential $\Psi_{2S_t}(x^t)$, and the attaining instance. We first comment on the statement.

The admissibility condition \eqref{eq:admissible} is the condition $\theta+L\gamma\le2$ of Lemmas~\ref{lemdet} and~\ref{lem:local} at every iteration, with $\gamma=\gamma_t$ and $\tau=2S_t$. It involves the future schedule, not future observations, and it permits early steps larger than $\frac{3}{2L}$ when a long tail of steps remains: at iteration $t$, it reads $L\gamma_t\le2-\frac{\gamma_t}{2S_t}$, and the subtracted term is small when $S_t\gg\gamma_t$. 
The last step, for which $S_{T-1}=\gamma_{T-1}$, necessarily satisfies $L\gamma_{T-1}\le\frac32=\beta_1$, the threshold of Section~\ref{sec:classical}.
No monotonicity of the steps is assumed. Adaptive stepsizes, chosen as functions of the observed iterates, are not covered: the parameter $2S_t$ of the potential would then depend on future observations, the vector $B(x^t)$ of Lemma~\ref{lem:local} would no longer be $\cF^t$-measurable, and the linear term in the noise would not vanish in conditional expectation. The upper bound \eqref{eq:upper} holds in every real Hilbert space; the condition $\dim\cH\ge T+1$ is needed only to build the attaining instance, which uses one coordinate for the initial distance and one coordinate per noise injection.
The lower bound holds for every positive schedule, admissible or not, and its instance attains the value $F_{L,T}$ in every case; without admissibility, however, this value need not be the worst-case value: for $T=1$, $\sigma_0=0$, and a step such that $2L\gamma(1-L\gamma)^2>1$, the quadratic $\ff=\frac L2\|\cdot\|^2$ with $\RR=0$ gives $\Delta_1=\frac{LD^2}{2}(1-L\gamma)^2>\frac{D^2}{4\gamma}$.\medskip

\noindent\textbf{Integrability.\ }
We first justify the expectations.
Since $-\nabla \ff(x^\star)\in\partial \RR(x^\star)$, the nonexpansiveness of the prox, the $\cF^t$-measurability of $x^t-x^\star-\gamma_t(\nabla \ff(x^t)-\nabla \ff(x^\star))$, and \eqref{eq:noise} give
\begin{align}
 \E\big[\norm{x^{t+1}-x^\star}^2\big]
 &\overset{\eqref{eq:algorithm},\eqref{eq:noise}}{\leq}
 \E\big[\norm{x^t-x^\star-\gamma_t(\nabla \ff(x^t)-\nabla \ff(x^\star))}^2\big]+\gamma_t^2\sigma_t^2\notag\\
 &\overset{\eqref{eq:interp}}{\leq}\E\big[\norm{x^t-x^\star}^2\big]+\gamma_t^2\sigma_t^2,\qquad \gamma_t\leq2/L.\label{eq:integrability}
\end{align}
The last inequality follows by expanding the square and using the cocoercivity of $\nabla\ff$, which is \eqref{eq:interp} added to the same inequality with $y$ and $z$ exchanged. This proves inductively that all iterates have finite second moments.
 
The mapping $B$ of Lemma~\ref{lem:local} is Lipschitz continuous, so that $B(x)$ has at most linear growth.
Therefore, its inner product with the current noise is integrable. The envelopes with a positive parameter are integrable by \eqref{eq:initial}; the pointwise bound \eqref{eq:local} also gives the integrability of the final objective gap, which is nonnegative.

\begin{proof}[Proof of the upper-bound part of Theorem~\ref{thm:exact}]
Set $\tau_t\eqdef2S_t$. Then $\tau_t\geq2\gamma_t$ and $\tau_t-2\gamma_t=\tau_{t+1}$. The local condition is exactly \eqref{eq:admissible}.
In particular, it implies $\gamma_t<2/L$, so that \eqref{eq:integrability} applies.
Moreover, with $\gamma=\gamma_t$ and $\tau=\tau_t$, the denominator $d$ of Lemma~\ref{lem:local} is $d_t$ in \eqref{eq:denominators}.
Apply the local inequality \eqref{eq:local} at $x=x^t$ with $e=e^t$.
Here $B(x^t)$ is $\cF^t$-measurable because the schedule is predetermined, so that the conditional expectation given $\cF^t$ of $\ip{B(x^t)}{e^t}$ vanishes;
then taking the total expectation and using \eqref{eq:noise} gives
\begin{equation}
 \E[\Psi_{\tau_{t+1}}(x^{t+1})-\Psi^\star]
 \overset{\substack{\eqref{eq:local},\eqref{eq:noise}\\\eqref{eq:c},\eqref{eq:denominators}}}{\leq}
 \E[\Psi_{\tau_t}(x^t)-\Psi^\star]
 +\frac{\sigma_t^2\gamma_t^2}{2d_t}.\label{eq:conditional}
\end{equation}
Summing yields
 $\displaystyle\E[\gap]
 \overset{\eqref{eq:conditional}}{\leq}
 \Psi_{2S_0}(x^0)-\Psi^\star
 +\frac12\sum_{t=0}^{T-1}\frac{\sigma_t^2\gamma_t^2}{d_t}\overset{\eqref{eq:initial},\eqref{eq:cost}}{\leq}F_{L,T}(\boldsymbol{\gamma};D,\boldsymbol{\sigma})$.
\end{proof}

\noindent\textbf{The attaining instance.\ }
It is built from the Huber function $\hub_\zeta$ of \eqref{eq:huber} and from absolute values, applied to orthonormal coordinates; a sum of Huber functions of orthonormal coordinates still has an $L$-Lipschitz gradient.

\begin{proof}[Proof of the attainment part of Theorem~\ref{thm:exact}]
Suppose first that $D>0$ and $\sigma_t>0$ for every $t$, and set, with
$\gc_t\eqdef\min\big\{\gamma_t,\frac1L\big\}$ and $d_t$ as in \eqref{eq:denominators},
\begin{equation}
 w_t\eqdef\frac{\gamma_t\sigma_t}{d_t},\qquad
 \zeta_t\eqdef\gc_tw_t,\qquad \lambda_t\eqdef(1-L\gc_t)w_t\geq0.\label{eq:lowerparameters}
\end{equation}
Let $E$ be a subspace of $\cH$ of dimension $T+1$, which exists since $\dim\cH\ge T+1$, equipped with an orthonormal basis, and let $z_b,z_0,\ldots,z_{T-1}$ denote the coordinates in this basis of the orthogonal projection of $x\in\cH$ onto $E$: one initial-distance coordinate and $T$ noise coordinates. The instance depends on these coordinates only, so that the iterates stay in $E$. Define
\begin{equation}
 \ff(z)\eqdef\sum_{t=0}^{T-1}\hub_{\zeta_t}(z_t),\qquad
 \RR(z)\eqdef\frac{D}{2S_0}|z_b|+\sum_{t=0}^{T-1}\lambda_t|z_t|,\label{eq:instance}
\end{equation}
and let $x^0$ be the point of $E$ with coordinates $z_b=D$ and $z_t=0$ for every $t$.
For a long step, that is, if $L\gamma_t>1$, we have $\lambda_t=0$: the coordinate $t$ then carries no nonsmooth term.
Every $\zeta_t$ is positive, so that the minimizers of $\Psi$ are the points whose $T+1$ coordinates vanish, even when some $\lambda_t$ are zero; $x^\star=0$ is one of them, with $\|x^0-x^\star\|=D$. At iteration $t$, the error
$e^t$ is the vector of $E$ whose only nonzero coordinate is $z_t=\pm\sigma_t$, the two signs being equiprobable and independent of the past. Both the unconditional bound in \eqref{eq:noise} and the stronger conditional variance bound hold with equality.

The initial-distance coordinate decreases by $\gamma_tD/(2S_0)$ at each step, stays positive, and ends at $D/2$. Its contribution is $D^2/(4S_0)$.
Before its injection, the noise coordinate $t$ is zero.
At injection, soft thresholding gives it the magnitude $\gamma_t(\sigma_t-\lambda_t)>0$. Thereafter, its sign is fixed,
and its magnitude decreases by $\gamma_jw_t$ at step $j>t$, since the composite slope on the linear Huber facet is $L\zeta_t+\lambda_t=w_t$. To verify that these updates remain on the linear Huber facet and do not cross the threshold, calculate the candidate final magnitude: 
since $\gamma_t\sigma_t=d_tw_t$, $d_t=2S_{t+1}+\omega(\gamma_t)=2S_{t+1}+\gc_t(2-L\gc_t)$, and
$(\gamma_t-\gc_t)(1-L\gc_t)=0$,
\begin{align}
 z_t^{\rm fin}
 &=\gamma_t(\sigma_t-\lambda_t)-S_{t+1}w_t\overset{\eqref{eq:lowerparameters}}{=}
 \big(d_t-\gamma_t(1-L\gc_t)-S_{t+1}\big)w_t
 \overset{\eqref{eq:denominators}}{=}(S_{t+1}+\gc_t)w_t
 \overset{\eqref{eq:lowerparameters}}{\geq}\zeta_t.\label{eq:finalcoordinate}
\end{align}
Every earlier candidate magnitude after injection is at least this value.
By induction, the claimed updates are therefore valid at every iteration. The contribution of this coordinate to the objective is
\begin{align}
 \hub_{\zeta_t}(z_t^{\rm fin})+\lambda_tz_t^{\rm fin}
 &\overset{\eqref{eq:huber},\eqref{eq:lowerparameters}}{=}
 w_tz_t^{\rm fin}-\frac L2\zeta_t^2\overset{\eqref{eq:finalcoordinate}}{=}
 \left(S_{t+1}+\gc_t-\frac L2\gc_t^2\right)w_t^2\overset{\eqref{eq:denominators},\eqref{eq:lowerparameters}}{=}
 \frac{\sigma_t^2\gamma_t^2}{2d_t}.\label{eq:coordinatevalue}
\end{align}
Summing these contributions proves the pathwise attainment of \eqref{eq:cost}.

If $\sigma_t=0$, replace the smooth term of that coordinate by $Lz_t^2/2$ and its nonsmooth term by zero, and inject no noise into it: the coordinate remains zero, and its term has a unique minimizer.
If $D=0$, use the same quadratic on the unused initial-distance coordinate and start it at zero. These changes preserve the value. If a unique minimizer is desired, add to $\ff$ the function $x\mapsto\frac L2\operatorname{dist}(x,E)^2$, where $\operatorname{dist}(x,E)$ is the distance from $x$ to $E$: $\nabla\ff$ remains $L$-Lipschitz, the iterates are unchanged, since they stay in $E$, and $0$ becomes the unique minimizer, in the degenerate cases as well.
\end{proof}

The instance \eqref{eq:instance} combines the two extremal mechanisms met in Section~\ref{sec:classical}. Its initial-distance coordinate is the instance of \citet[Theorem~4.1]{thg17}, $\ff=0$ and $\RR=\lambda|\cdot|$ in dimension one, which attains the deterministic rate $\frac{D^2}{4S_0}$. 
When $L\gamma_t\ge1$, the coefficient $\lambda_t$ vanishes, so that the noise coordinate is a Huber coordinate, as in the instance of \citet{gs26}; when $L\gamma_t<1$, the noise coordinate combines curvature and a kink, in the proportions $L\gamma_t$ and $1-L\gamma_t$ of the composite slope $w_t=L\zeta_t+\lambda_t$ of \eqref{eq:lowerparameters}, as in the one-step instance discussed at the end of Section~\ref{sec:noisystep}.
The exactness \eqref{eq:exact} thus holds on the whole admissible range, which contains the range $\gamma_t\le\frac{3}{2L}$ of the smooth result of \citet{gs26}.

\medskip\noindent\textbf{A two-step example.\ }
Let us follow the instance in the smallest nontrivial case: $T=2$, $L=D=\sigma_0=\sigma_1=1$, and $\gamma_0=\gamma_1=\frac12$. Then $S_0=1$, $S_1=\frac12$, and $\omega\big(\frac12\big)=\frac34$, so that $d_0=\frac74$ and $d_1=\frac34$, and the cost is
\begin{equation}
 F_{1,2}\overset{\eqref{eq:cost}}{=}\frac14+\frac1{14}+\frac16=\frac{41}{84}.\label{eq:workedcost}
\end{equation}
The two injections have the same variance bound, but the first one costs $\frac1{14}$ and the last one $\frac16$, because no iteration remains to damp it. In the instance, $w_0=\frac27$ and $w_1=\frac23$, so that $\zeta_0=\lambda_0=\frac17$ and $\zeta_1=\lambda_1=\frac13$, by \eqref{eq:lowerparameters}. Starting from $(z_b,z_0,z_1)=(1,0,0)$, the magnitudes of the three coordinates are, whatever the signs of the noise,
\begin{equation*}
 (1,0,0)\ \overset{\eqref{eq:algorithm},\eqref{eq:lowerparameters}}{\longmapsto}\ \big(\tfrac34,\tfrac37,0\big)\ \overset{\eqref{eq:algorithm},\eqref{eq:lowerparameters}}{\longmapsto}\ \big(\tfrac12,\tfrac27,\tfrac13\big).
\end{equation*}
Indeed, the initial-distance coordinate decreases by $\gamma_t\frac{D}{2S_0}=\frac14$ at each step; at the first step, soft thresholding gives $|z_0|=\frac12\big(1-\frac17\big)=\frac37$; the second step subtracts $\frac12w_0=\frac17$ from it and gives $|z_1|=\frac12\big(1-\frac13\big)=\frac13$. At the end, the objective value is $\frac{D}{2S_0}\cdot\frac12=\frac14$ in the first coordinate and, by \eqref{eq:coordinatevalue}, $\frac1{14}$ and $\frac16$ in the other two: these are the three terms of \eqref{eq:workedcost}. Thus the formula describes an actual run, and not only a sum of estimates that are sharp separately.

\begin{corollary}[Constant steps and deterministic specialization]\label{cor:constant}
Let $\gamma_t\equiv\gamma\leq \frac3{2L}$ and $\sigma_t\equiv\sigma$. Then
\begin{align}
 \E[\gap]&\overset{\eqref{eq:upper},\eqref{eq:cost},\eqref{eq:denominators}}{\leq}\frac{D^2}{4\gamma T}
 +\frac{\gamma^2\sigma^2}{2}\sum_{k=0}^{T-1}\frac1{2k\gamma+\omega(\gamma)}\label{eq:constant}\\
 &\overset{\eqref{eq:constant}}{\leq}\frac{D^2}{4\gamma T}+\frac{\gamma^2\sigma^2}{2\omega(\gamma)}
 +\frac{\gamma\sigma^2}{4}\log\Big(1+\frac{2(T-1)\gamma}{\omega(\gamma)}\Big).\notag
\end{align}
The first bound is the exact worst-case value if $\dim\cH\ge T+1$.
If $\sigma_t=0$ for every $t$, the bound $D^2/(4S_0)$ is exact for every schedule satisfying \eqref{eq:admissible}, including constant stepsizes up to $\frac{3}{2L}$.
\end{corollary}
\begin{proof}
Substitution into \eqref{eq:upper} gives the first bound, since $d_t=2(T-t-1)\gamma+\omega(\gamma)$ by \eqref{eq:denominators}; the index $k=T-t-1$ is the number of iterations that remain after the injection. The second bound follows by bounding the decreasing summands after the first one by an integral. The exactness is \eqref{eq:exact}.
In the deterministic case, the one-dimensional instance $\ff=0$, $\RR(z)=D|z|/(2S_0)$, $x^0=D$ gives $x^T=D/2$ for every positive schedule, and hence attains the bound.
If $D=0$, the zero trajectory suffices.
\end{proof}

If $L\gamma\le1$, then $\omega(\gamma)=\gamma(2-L\gamma)$ and the noise term of \eqref{eq:constant} is
$\frac{\gamma\sigma^2}{2}\sum_{k=0}^{T-1}\frac{1}{2k+2-L\gamma}$;
if $1\le L\gamma\le\frac32$, then $\omega(\gamma)=\frac1L$ and it is
$\frac{L\gamma^2\sigma^2}{2}\sum_{k=0}^{T-1}\frac{1}{1+2kL\gamma}$, the noise term of \eqref{eq:gs}. For a fixed stepsize $\gamma\le\frac{3}{2L}$ and $\sigma>0$, the exact worst-case value \eqref{eq:constant} grows logarithmically with $T$.
This is a worst-case statement, with an instance that depends on the horizon, since the instance of Theorem~\ref{thm:exact} has $T+1$ coordinates; it is not the divergence of a fixed trajectory, nor of the iterates on a fixed problem. Section~\ref{sec:design} shows how the choice of the schedule changes this conclusion.

\begin{remark}[Initial-error term]\label{rembias}
The term $\frac{D^2}{4\gamma T}$ is exact: as recalled in Section~\ref{sec:classical}, it is attained in the case $\ff=0$ with $\RR=\lambda|\cdot|$ in dimension one. It was known to be the exact worst case of \algn{Prox-GD} for $\gamma\le\frac1L$ \citep[Theorems 7 and 8]{tv23}; Corollary~\ref{cor:constant} with $\sigma=0$ extends this to $\gamma\le\frac{3}{2L}$ and, more generally, shows that $\frac{D^2}{4S_0}$ is the exact worst case of \algn{Prox-GD} for every admissible schedule. Compared with the smooth initial-error term $\frac{D^2}{4\gamma T+2/L}$ of \eqref{eq:gs}, the composite class loses the $2/L$: the envelope of $\Psi$ is only $\frac{1}{\tau}$-smooth, so that \eqref{eq:initial} cannot see $L$ (Section~\ref{sec:envelopefacts}).
Compared with the classical bound \eqref{eqr1}, the constant is improved by a factor of 2,
and the stepsize range is extended from $\frac1L$ to $\frac{3}{2L}$.
\end{remark}

\begin{remark}[Noise term]\label{remnoise}
The noise term has the same structure as in the smooth case \eqref{eq:gs}: a logarithm in $T$ with the constant $\frac{\gamma\sigma^2}{4}$ in front,
plus an $O(\gamma\sigma^2)$ term accounting for the last injections.
For $1\le L\gamma\le\frac32$, it \emph{is} the smooth noise term: \eqref{eq:constant} and \eqref{eq:gs} differ by their first terms only. For $L\gamma<1$, an injection followed by $k$ more iterations has the weight $\frac{1}{2k+2-L\gamma}$, in place of the smaller smooth weight $\frac{L\gamma}{1+2kL\gamma}=\frac{1}{2k+1/(L\gamma)}$; for small $L\gamma$ and $TL\gamma\gg1$,
the smooth sum is $\approx\frac12\log(1+2TL\gamma)$, while the composite one is $\approx\frac12(\log T+C_{\rm E})$,
where $C_{\rm E}\approx0.577$ is Euler's constant, a gap of about $\frac{\gamma\sigma^2}{4}\big(\log\frac{1}{2L\gamma}+C_{\rm E}\big)$, which is the price of the prox discussed in Section~\ref{sec:noisystep}; the identity \eqref{eq:price} below gives this price exactly, for every schedule. By Theorem~\ref{thm:exact}, this price is real and not an artifact of the analysis: the composite weights are exact.
In the case $\ff=0$, which is $L$-smooth for every $L>0$, so that we can let $L\to0$ in \eqref{eq:constant} and \eqref{eq:exact}, we obtain the exact bound $\E[\Delta_T]\le\frac{D^2}{4\gamma T}+\frac{\gamma\sigma^2}{4}\sum_{t=1}^T\frac1t$,
attained when $\RR$ is a separable $\ell_1$ function with weights $\frac{\sigma}{2t}$,
$t=1,\ldots,T$; see Corollary~\ref{cor:ppa} below for general schedules. The bound \eqref{eq:lz} of \citet[Lemma 4.2]{lz24}, $\frac{D^2}{2\gamma T}+2\gamma\sigma^2\sum_{t=1}^T\frac1t$ for $\gamma\le\frac{1}{2L}$,
is improved by a factor of 2 on the initial-error term and by a factor between 6 and 8 on the noise term;
moreover, our bounds hold for stepsizes up to $\frac{3}{2L}$.
This is a comparison of the two bounds for the same constant stepsize, in the range $\gamma\le\frac{1}{2L}$ where both apply; it says nothing about the results of \citet{lz24} for decaying stepsizes or in high probability, which have no counterpart here.
\end{remark}

\noindent\textbf{The exact composite penalty.\ }
The composite cost \eqref{eq:cost} and the smooth cost \eqref{eq:smooth} of Proposition~\ref{prop:smooth} can be compared term by term, for every positive schedule. 
Set $\tilde d_t\eqdef2S_{t+1}+\frac1L$, so that $\tilde d_t-d_t=\frac1L-\omega(\gamma_t)=\frac1L(1-L\gamma_t)_+^2$ by \eqref{eq:omega};
this difference vanishes for a long step. A direct subtraction gives the exact identity
\begin{align}
 F_{L,T}-F_{L,T}^{\mathrm{smooth}}
 &\overset{\eqref{eq:cost},\eqref{eq:smooth},\eqref{eq:omega}}{=}
 \frac{D^2/L}{2S_0(4S_0+2/L)}
 +\frac1{2L}\sum_{t=0}^{T-1}
 \frac{\sigma_t^2\gamma_t^2(1-L\gamma_t)_+^2}{d_t\tilde d_t}.\label{eq:price}
\end{align}
The first term is the loss of the $2/L$ in the initial-error term. The variance coefficients agree exactly at the iterations where $L\gamma_t\ge1$, and the composite one is larger at the short steps, as for the one-step constants compared at the end of Section~\ref{sec:noisystep}. For an admissible schedule, both costs are exact worst-case values, so that \eqref{eq:price} is the exact price of the prox, and not the difference of two proof bounds: for the composite cost, this is Theorem~\ref{thm:exact}; for the smooth cost, the upper bound is Proposition~\ref{prop:smooth},
and the lower bound is that of \citet{gs26} for a uniform variance bound,
whose instance adapts to a nonuniform profile, when $\dim\cH\ge T+1$: in the noise coordinate $t$, use the Huber function $\hub_{\zeta_t}$ with $\zeta_t\eqdef\gamma_t\sigma_t/(L\tilde d_t)$
and inject $\pm\sigma_t$ at iteration $t$,
so that the coordinate stays on the linear facet after the injection, ends at magnitude $\big(\frac1L+S_{t+1}\big)L\zeta_t\geq\zeta_t$, and contributes $\frac{\sigma_t^2\gamma_t^2}{2\tilde d_t}$; in the initial-distance coordinate, start at $D$ and use the threshold $\frac{D}{1+2LS_0}$, which contributes $\frac{D^2}{4S_0+2/L}$ (the coordinates with $\sigma_t=0$, and the initial-distance coordinate if $D=0$, are omitted).

The case $\ff=0$ is included in Assumption~\ref{ass:model} for every $L>0$, and the cost \eqref{eq:cost} decreases when $L$ decreases, so that its exact performance is obtained in the limit $L\to0$. It is also the case in which the admissibility condition disappears: \eqref{eq:admissible} becomes $\frac{\gamma_t}{2S_t}\le2$, which always holds.
\begin{corollary}[Exact performance when $\ff=0$]\label{cor:ppa}
Let $\ff=0$, so that \eqref{eq:algorithm} reads
$x^{t+1}=\prox_{\gamma_t \RR}(x^t-\gamma_te^t)$,
let $\RR$ be proper, lower semicontinuous, and convex,
with a minimizer $x^\star$ such that $\norm{x^0-x^\star}\le D$, and let the errors satisfy \eqref{eq:noise}. Let $\boldsymbol{\gamma}$ be any positive predetermined schedule. Then $\E[\gap]\le F_{0,T}(\boldsymbol{\gamma};D,\boldsymbol{\sigma})$, where
\begin{equation}
 F_{0,T}(\boldsymbol{\gamma};D,\boldsymbol{\sigma})\eqdef\frac{D^2}{4S_0}+\frac14\sum_{t=0}^{T-1}\frac{\sigma_t^2\gamma_t^2}{S_t},\label{eq:ppacost}
\end{equation}
and this bound is the exact worst-case value if $\dim\cH\ge T+1$.
\end{corollary}
\begin{proof}
For the upper bound, note that $\ff=0$ is $\varepsilon$-smooth for every $\varepsilon>0$.
For a fixed finite schedule, all steps satisfy \eqref{eq:admissible} for sufficiently small $\varepsilon$. Apply Theorem~\ref{thm:exact} and let $\varepsilon\downarrow0$;
the algorithm and the left-hand side of \eqref{eq:upper} do not depend on $\varepsilon$.
For attainment, use $\RR(z)=D|z_b|/(2S_0)+\sum_t\lambda_t|z_t|$ with $\lambda_t=\gamma_t\sigma_t/(2S_t)$ and the same independent one-coordinate injections. The calculation in \eqref{eq:finalcoordinate}--\eqref{eq:coordinatevalue} applies with $L=0$. If a unique minimizer is desired, add to $\RR$ a positive quadratic function of every coordinate $z_t$ with $\sigma_t=0$, and of $z_b$ if $D=0$, as well as the function $x\mapsto\alpha\operatorname{dist}(x,E)^2$, for some $\alpha>0$, where $E$ is the subspace of the construction: the iterates stay in $E$, these coordinates remain zero, the objective values are unchanged, and $0$ becomes the unique minimizer, also when $\cH$ is larger than $E$.
\end{proof}
The cost $F_{0,T}$ is the limit $L\to0$ of $F_{L,T}$, and 
the case $\ff=0$ is the only one in which the exact worst case is known for every positive schedule, without restriction on the stepsizes.
It serves as a reference in the rest of the paper: Section~\ref{sec:design} computes its optimal schedule in closed form, and Section~\ref{sec:asymptotic} shows that this schedule achieves the asymptotic minimax law of Theorem~\ref{thm:secondorder} for every $L>0$.

\section{Exact finite-horizon schedule design}\label{sec:design}
Since the cost of a schedule is known exactly, the schedule can be optimized. Throughout this section, $D,\sigma>0$ and the variance bound is uniform, $\sigma_t\equiv\sigma$, except in Section~\ref{sec:nonuniform}. The positivity of $D$ and $\sigma$ is essential: if $D=0<\sigma$, letting the total stepsize vanish drives the cost to zero, and if $\sigma=0<D$, letting the total stepsize grow does the same, so that neither problem has a minimizer; if $D=\sigma=0$, every schedule has zero cost. 
Throughout this section and the next one, we suppose that $\dim\cH\ge T+1$, for every horizon $T$ considered, so that, by Theorem~\ref{thm:exact},
$F_{L,T}$ is a lower bound on the worst case of every positive schedule and the exact worst-case value of every admissible one,
and $F_{0,T}$ is the exact worst-case value of every positive schedule in the case $\ff=0$, by Corollary~\ref{cor:ppa}. The design problem is to find a schedule with the smallest worst case, that is, to determine the minimax value
\begin{equation}
 M_{L,T}\eqdef \inf_{\boldsymbol{\gamma}\in(0,\infty)^T}\cW_{L,T}(\boldsymbol{\gamma};D,\sigma\mathbf1)
 \label{eq:globalvalue}
\end{equation}
and a schedule that attains it. Sections~\ref{sec:convex}--\ref{sec:scalar} 
study the minimization of the cost $F_{L,T}$ over all positive schedules:
the problem is strictly convex, its minimizer is strictly decreasing, and it can be computed by a scalar search. Section~\ref{sec:global} shows that this minimizer solves the minimax problem whenever it is admissible, which is the case unless the noise level is low.
Section~\ref{sec:ppa} solves the case $\ff=0$ in closed form, and Section~\ref{sec:nonuniform} extends this solution to a nonuniform variance profile.

\subsection{A strictly convex problem}\label{sec:convex}
Define the optimal value of the cost over all positive schedules,
\begin{equation}
 V_{L,T}\eqdef \inf_{\boldsymbol{\gamma}\in(0,\infty)^T}F_{L,T}(\boldsymbol{\gamma};D,\sigma\mathbf1).
 \label{eq:optvalue}
\end{equation}
No upper bound is imposed on the stepsizes; in particular, the schedules are not required to be admissible.

\begin{theorem}[Existence, uniqueness, and optimality conditions]\label{thm:convex}
The cost $F_{L,T}(\,\cdot\,;D,\sigma\mathbf1)$ is a strictly convex and continuously differentiable function on $(0,\infty)^T$, and the infimum in \eqref{eq:optvalue} is attained at a unique positive schedule $\boldsymbol{\gamma}^\star$. With $d_i$ as in \eqref{eq:denominators}, the partial derivative of the cost with respect to $\gamma_k$, denoted by $\partial_kF$, is
\begin{equation}
 \partial_kF\overset{\eqref{eq:cost},\eqref{eq:denominators}}{=}-\frac{D^2}{4S_0^2}
 +\sigma^2\left\{\frac{\gamma_k\big(2S_{k+1}+\min\{\gamma_k,\frac1L\}\big)}{d_k^2}
                 -\sum_{i<k}\frac{\gamma_i^2}{d_i^2}\right\},
 \label{eq:gradient}
\end{equation}
so that a positive schedule is the minimizer $\boldsymbol{\gamma}^\star$ if and only if
\begin{equation}
 \partial_kF=0,\quad\forall k=0,\ldots,T-1.\label{eq:kkt}
\end{equation}
\end{theorem}
The proof is given in Appendix~\ref{app:proofconvex}; its main points are the following. For a step $\gamma_i\le\frac1L$, the noise term is an increasing convex function of a quadratic-over-linear function of the schedule:
\begin{equation}
 r_i\eqdef\frac{\gamma_i^2}{2S_i},\qquad
 \phi(r)\eqdef\frac r{1-Lr},\qquad
 \frac{\gamma_i^2}{d_i}\overset{\eqref{eq:denominators}}{=}\frac{\gamma_i^2}{2S_i-L\gamma_i^2}
 =\phi(r_i).\label{eq:convexcomposition}
\end{equation}
For a step $\gamma_i\ge\frac1L$, the noise term $\gamma_i^2/\big(2S_{i+1}+\frac1L\big)$ is itself quadratic-over-linear. These two convex pieces match in a continuously differentiable way along $\{L\gamma_i=1\}$, because $\omega'\big(\frac1L\big)=0$, and this preserves convexity. A minimizer exists because the initial-error term prevents the total stepsize from vanishing and the noise terms prevent the steps from growing; a first-order argument shows that none of its steps is zero.

In the sequel, we call $\boldsymbol{\gamma}^\star$ the \emph{optimal schedule}, for short: it is the schedule with the smallest cost, and it is minimax optimal, that is, it has the smallest worst case, whenever it is admissible (Theorem~\ref{thm:global}).

The conditions \eqref{eq:kkt} characterize the minimizer. When a schedule satisfies them only approximately, for instance because it has been computed numerically, the following proposition bounds its suboptimality.

\begin{proposition}[A computable bound on the optimization error]\label{prop:numericalcertificate}
Let $L,D,\sigma>0$, let $\bar{\boldsymbol{\gamma}}$ be a positive schedule, with cost $\bar F\eqdef F_{L,T}(\bar{\boldsymbol{\gamma}};D,\sigma\mathbf1)$, and let $\partial_tF$ denote the partial derivatives \eqref{eq:gradient} of the cost at $\bar{\boldsymbol{\gamma}}$. Let $U_{T-1},\ldots,U_0$ be positive numbers chosen backward, such that
\begin{equation}
 U_t\geq\sqrt{\frac{2\bar F}{\sigma^2}\Big(2\sum_{j>t}U_j+\frac1L\Big)},\quad\forall t=T-1,\ldots,0.
 \label{eq:sublevelbox}
\end{equation}
Then
\begin{equation}
 0\leq \bar F-V_{L,T}\leq
 \mathcal E(\bar{\boldsymbol{\gamma}})\eqdef
 \sum_{t=0}^{T-1}\Big[(\partial_tF)_+\,\bar\gamma_t+(-\partial_tF)_+\,(U_t-\bar\gamma_t)\Big].
 \label{eq:numericalgap}
\end{equation}
Consequently, for every admissible schedule $\boldsymbol{\eta}$,
\begin{equation}
 \bar F-\mathcal E(\bar{\boldsymbol{\gamma}})
 \overset{\eqref{eq:numericalgap},\eqref{eq:universallower}}{\leq} M_{L,T}\overset{\eqref{eq:exact},\eqref{eq:globalvalue}}{\leq} F_{L,T}(\boldsymbol{\eta};D,\sigma\mathbf1).
 \label{eq:certifiedminimax}
\end{equation}
In particular, if $\bar{\boldsymbol{\gamma}}$ is admissible, its worst case exceeds the minimax value $M_{L,T}$ by at most $\mathcal E(\bar{\boldsymbol{\gamma}})$.
\end{proposition}
\begin{proof}
Since $\omega\le\frac1L$, every positive schedule $\boldsymbol{\gamma}$ whose cost is at most $\bar F$ satisfies
\begin{equation}
 \gamma_t^2\overset{\eqref{eq:cost},\eqref{eq:denominators},\eqref{eq:omega}}{\leq}
 \frac{2\bar F}{\sigma^2}\Big(2S_{t+1}+\frac1L\Big),\quad\forall t=0,\ldots,T-1.
 \label{eq:coordinatecap}
\end{equation}
By backward induction, \eqref{eq:coordinatecap} and \eqref{eq:sublevelbox} give $0<\gamma_t\le U_t$ for every $t$. This applies to $\bar{\boldsymbol{\gamma}}$ and to the minimizer $\boldsymbol{\gamma}^\star$ of Theorem~\ref{thm:convex}. Since the cost is convex and differentiable,
\begin{equation*}
 \bar F-V_{L,T}\leq\sum_{t=0}^{T-1}\partial_tF\,(\bar\gamma_t-\gamma_t^\star)
 \overset{\eqref{eq:sublevelbox},\eqref{eq:coordinatecap}}{\leq}\sum_{t=0}^{T-1}\Big[(\partial_tF)_+\,\bar\gamma_t+(-\partial_tF)_+\,(U_t-\bar\gamma_t)\Big],
\end{equation*}
where the second inequality holds because $0<\gamma_t^\star\le U_t$. This is \eqref{eq:numericalgap}. By Theorem~\ref{thm:exact}, the worst case of every positive schedule is at least its cost, so that $M_{L,T}\ge V_{L,T}$, and the worst case of an admissible schedule $\boldsymbol{\eta}$ is its cost, so that $M_{L,T}\le F_{L,T}(\boldsymbol{\eta};D,\sigma\mathbf1)$. With \eqref{eq:numericalgap}, this proves \eqref{eq:certifiedminimax}, and the last assertion is the case $\boldsymbol{\eta}=\bar{\boldsymbol{\gamma}}$.
\end{proof}
Proposition~\ref{prop:numericalcertificate} converts the gradient of the cost at a candidate schedule into a bound on its suboptimality. A small gradient alone would not provide such a bound; the numbers $U_t$ supply the missing information, a box that contains the minimizer.

\subsection{The optimal steps decrease}\label{sec:order}
The exact formula gives more than the first-order conditions \eqref{eq:kkt}: the exchange of two adjacent steps reveals the order of the optimal steps without solving these conditions.

\begin{theorem}[Order of the optimal steps]\label{thm:order}
The minimizer $\boldsymbol{\gamma}^\star$ of Theorem~\ref{thm:convex} is strictly decreasing:
\begin{equation}
 \gamma_0^\star>\gamma_1^\star>\cdots>\gamma_{T-1}^\star>0.\label{eq:ordershape}
\end{equation}
\end{theorem}
\begin{proof}
Consider adjacent steps $a,b$ with a fixed tail $Y\geq0$ after them, and put $S\eqdef a+b+Y$.
The earlier terms and the initial-error term in \eqref{eq:cost} are not changed by the exchange of $a$ and $b$.
Define, with $S,Y$ fixed and $\omega$ as in \eqref{eq:omega},
\begin{equation}
 \xi_{S,Y}(z)\eqdef \frac{z}{[2Y+\omega(z)][2S-2z+\omega(z)]}.
 \label{eq:exchangefunction}
\end{equation}
Direct subtraction gives
\begin{equation*}
 F_{L,T}(\ldots,a,b,\ldots)-F_{L,T}(\ldots,b,a,\ldots)
 \overset{\eqref{eq:cost},\eqref{eq:exchangefunction}}{=}
 \sigma^2ab\,[\xi_{S,Y}(b)-\xi_{S,Y}(a)].
\end{equation*}
On the interval between $a$ and $b$, we have $z<S$, so that the two factors of the denominator are positive. The function $z\mapsto z/[2Y+\omega(z)]$ is increasing: its derivative is $(2Y+Lz^2)/[2(z+Y)-Lz^2]^2>0$ if $Lz<1$ and $1/[2Y+\frac1L]>0$ if $Lz>1$. The function $z\mapsto[2S-2z+\omega(z)]^{-1}$ is positive and increasing as well, since $2S-2z+\omega(z)$ is equal to $2S-Lz^2$ if $Lz\le1$ and to $2S-2z+\frac1L$ if $Lz\ge1$. Hence $\xi_{S,Y}$ is strictly increasing.
An inversion $a<b$ would therefore strictly increase the cost, which proves that the minimizer is nonincreasing.

To exclude ties, suppose $\gamma_i=\gamma_{i+1}=a$, and put $Y\eqdef S_{i+2}$ and $S\eqdef2a+Y$. If $La\le1$, then $d_i=2S-La^2$ and $d_{i+1}=2(S-a)-La^2$, and subtraction of adjacent gradient components gives
\begin{align}
 \partial_iF-\partial_{i+1}F
 &\overset{\eqref{eq:gradient}}{=}
 \sigma^2a\left(\frac{2S}{d_i^2}-\frac{2S-3a}{d_{i+1}^2}\right)\notag\\
 &=-\frac{\sigma^2a^2[4SY+La^2(4S-3La^2)]}{d_i^2d_{i+1}^2}<0.
 \label{eq:tiederivative}
\end{align}
Indeed, $La\leq1$ and $S\geq2a$ imply $4S-3La^2>0$.
If $La\ge1$, then $d_{i+1}=2Y+\frac1L$ and $d_i=d_{i+1}+2a$, and
\begin{equation}
 \partial_iF-\partial_{i+1}F
 \overset{\eqref{eq:gradient}}{=}
 \sigma^2a\left(\frac{d_i+a}{d_i^2}-\frac{1}{d_{i+1}}\right)
 =-\frac{\sigma^2a^2(d_{i+1}+4a)}{d_i^2d_{i+1}}<0.
 \label{eq:tiederivativelong}
\end{equation}
In both cases, this contradicts $\partial_iF=\partial_{i+1}F=0$, which holds at the minimizer by \eqref{eq:kkt}. This proves \eqref{eq:ordershape}.
\end{proof}
Theorem~\ref{thm:order} proves the monotonicity of the optimal schedule at every positive noise level. For the steps larger than $\frac1L$, the stationarity conditions take a simple form: if $\gamma_t\ge\frac1L$ and $\gamma_{t+1}\ge\frac1L$, the equality $\partial_tF=\partial_{t+1}F$ reads $\frac{\gamma_{t+1}}{d_{t+1}}=\frac{\gamma_t}{d_t}+\big(\frac{\gamma_t}{d_t}\big)^2$. This is the recurrence that defines the minimax optimal schedule of the smooth case \citep{gs26}, as it should be, since the cost of such steps is the smooth one.
The strict decrease of the last steps contrasts with the case $\ff=0$, where the last two optimal steps are equal (Theorem~\ref{thm:ppaopt} below).

\subsection{Constructing the optimizer by one scalar search}\label{sec:scalar}

The ordering theorem makes a direct construction possible: once the last step is fixed, the equality of the gradient components at adjacent steps determines every earlier step. The last step is then found by a scalar search.

\begin{theorem}[A scalar construction of the optimizer]\label{thm:scalarconstruction}
Let $L,D,\sigma>0$. For each $s>0$, construct a schedule $\Gamma(s)$ backward. Set $\Gamma_{T-1}(s)\eqdef s$. Given $b=\Gamma_{i+1}(s)$ and $Y=\sum_{j>i}\Gamma_j(s)$, let $\Gamma_i(s)$ be the unique $a>b$ satisfying $P_Y(a)=B_Y(b)$, where, with $\omega$ as in \eqref{eq:omega},
\begin{equation}
 P_Y(a)\eqdef \frac{a\big(a+\min\{a,\frac1L\}+2Y\big)}{[2Y+\omega(a)]^2},\qquad
 B_Y(b)\eqdef \frac{b\big(2Y-2b+\min\{b,\frac1L\}\big)}{[2Y-2b+\omega(b)]^2}.
 \label{eq:scalarfunctions}
\end{equation}
Finally, put
\begin{equation*}
 \Lambda(s)\eqdef \partial_{T-1}F(\Gamma(s)),
\end{equation*}
where the gradient is \eqref{eq:gradient}.
Then $\Lambda$ has a unique zero $s^\star>0$; it is negative below this zero and positive above it, and $\Gamma(s^\star)$ is the minimizer $\boldsymbol{\gamma}^\star$ of Theorem~\ref{thm:convex}. Bisection finds this zero; the monotonicity of $\Lambda$ is neither needed nor asserted.
\end{theorem}
The functions \eqref{eq:scalarfunctions} are such that, for every positive schedule and every pair of adjacent steps, with $Y=S_{i+1}$,
\begin{equation}
 \partial_iF-\partial_{i+1}F\overset{\eqref{eq:gradient},\eqref{eq:scalarfunctions}}{=}
 \sigma^2[P_Y(\gamma_i)-B_Y(\gamma_{i+1})].
 \label{eq:adjacentstationarity}
\end{equation}
The construction thus makes all the gradient components of $\Gamma(s)$ equal to $\Lambda(s)$,
and a zero of $\Lambda$ gives a schedule that satisfies \eqref{eq:kkt}. The existence and uniqueness of every step and of the zero are proved in Appendix~\ref{app:proofscalar}.

Figure~\ref{fig:schedules} shows the minimizer $\boldsymbol{\gamma}^\star$ for $L=D=1$ and $T=30$. For $\sigma=0.2$ and $\sigma=1$, all its steps are smaller than $\frac1L$. For $\sigma=0.1$, its first $14$ steps are larger than $\frac1L$, the first one being $\frac{1.70}{L}$; this schedule is admissible, and its cost $0.01731$ is its exact worst case, whereas the best schedule with steps at most $\frac1L$ costs $0.01775$.
By Theorem~\ref{thm:global} below, these three schedules are the minimax optimal schedules among all positive predetermined schedules. The right panel compares the optimizer, at $\sigma=0.2$, with the optimal schedule of the case $\ff=0$ (Theorem~\ref{thm:ppaopt} below) and with the linearly decaying schedule \eqref{eq:linearschedule} of Section~\ref{sec:asymptotic}: the first is almost indistinguishable from the optimizer, its cost being $0.03418$ instead of $0.03416$, whereas linear decay, which ends with the two steps $2a,a$, costs $0.03456$.

The nearly linear decay of the optimal steps in Figure~\ref{fig:schedules} has a simple explanation. In the case $\ff=0$, the noise term of the cost \eqref{eq:ppacost} is $\frac{\sigma^2}{4}\sum_t\gamma_t^2/S_t$, and its continuous-time counterpart, with $\gamma=-S'$, is
\begin{equation*}
 \frac{\sigma^2}{4}\int_0^T\frac{S'(t)^2}{S(t)}\,dt=\sigma^2\int_0^T\Big(\frac{d}{dt}\sqrt{S(t)}\Big)^2dt.
\end{equation*}
For a given total stepsize $S(0)$, with $S(T)=0$, this energy is minimized when $\sqrt S$ is an affine function of $t$, that is, when the stepsize decays linearly to zero; and for $L>0$, the denominators $d_t$ are close to $2S_t$. The discrete optimal schedule departs from this picture mainly at its last step, which is damped by itself only and is therefore larger than linear decay would make it, as can be seen at $t=T-1$ in the left panel.

\begin{figure}[t]\centering
\includegraphics[width=\linewidth]{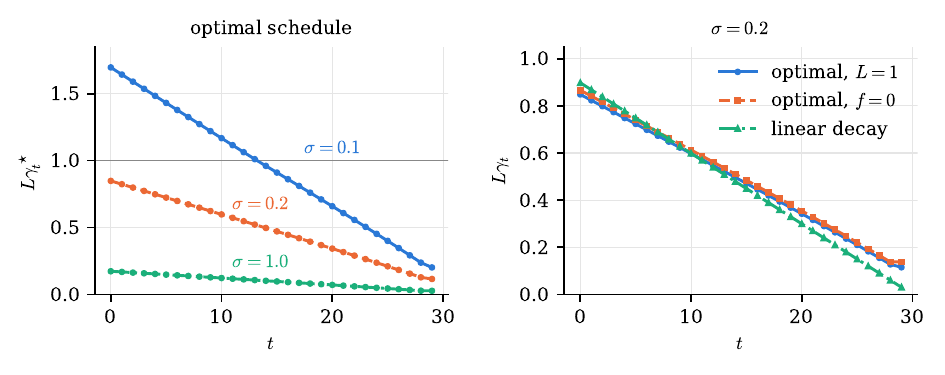}
\caption{Left: the optimal schedule $\boldsymbol{\gamma}^\star$, minimizer of the cost over all positive schedules (Theorems~\ref{thm:order} and~\ref{thm:scalarconstruction}), for $L=D=1$, $T=30$, and three noise levels $\sigma$: 
the steps are strictly decreasing, and they can be larger than $\frac1L$, the level of the horizontal line. The three schedules are admissible.
Right: for $\sigma=0.2$, the same optimizer, the optimal schedule $\widehat{\boldsymbol{\gamma}}$ of the case $\ff=0$ (Theorem~\ref{thm:ppaopt}), and the linearly decaying schedule \eqref{eq:linearschedule}.}\label{fig:schedules}
\end{figure}

\subsection{Minimax optimal schedules}\label{sec:global}
The cost $F_{L,T}$ is a lower bound on the worst case of every positive schedule,
and it is the worst case of the admissible ones (Theorem~\ref{thm:exact}). Its minimizer over all positive schedules therefore solves the minimax problem \eqref{eq:globalvalue} whenever it is admissible.

\begin{theorem}[Minimax optimal schedules among all predetermined schedules]\label{thm:global}
Let $L,D,\sigma>0$ and $\sigma_t\equiv\sigma$, and let $\boldsymbol{\gamma}^\star$ be the minimizer of $F_{L,T}(\,\cdot\,;D,\sigma\mathbf1)$ over all positive schedules, 
given by Theorems~\ref{thm:convex}--\ref{thm:scalarconstruction}.
\begin{enumerate}[label=(\roman*)]
\item $V_{L,T}\le M_{L,T}\le F_{L,T}(\boldsymbol{\gamma};D,\sigma\mathbf1)$ for every admissible schedule $\boldsymbol{\gamma}$.
\item If $\boldsymbol{\gamma}^\star$ is admissible, then
\begin{equation}
 M_{L,T}=V_{L,T}=F_{L,T}(\boldsymbol{\gamma}^\star;D,\sigma\mathbf1),\label{eq:globalexact}
\end{equation}
and $\boldsymbol{\gamma}^\star$ is the unique positive schedule whose worst case is $M_{L,T}$.
\item The largest step of $\boldsymbol{\gamma}^\star$ satisfies
\begin{equation}
 L\gamma_0^\star\le\max\bigg\{1,\ \frac14\bigg(1+\sqrt{1+\frac{2L^2D^2}{\sigma^2}}\,\bigg)\bigg\}.\label{eq:firststepbound}
\end{equation}
Consequently, all the steps of $\boldsymbol{\gamma}^\star$ are at most $\frac{3}{2L}$, so that it is admissible, if $\sigma\ge\frac{LD}{2\sqrt3}$.
\end{enumerate}
\end{theorem}
\begin{proof}
(i) By Theorem~\ref{thm:exact}, $\cW_{L,T}(\boldsymbol{\gamma};D,\sigma\mathbf1)\ge F_{L,T}(\boldsymbol{\gamma};D,\sigma\mathbf1)$ for every positive schedule, and taking the infimum gives the first inequality. 
The worst case of an admissible schedule is its cost, by \eqref{eq:exact}, and $M_{L,T}$ is at most this worst case; this gives the second inequality.

(ii) If $\boldsymbol{\gamma}^\star$ is admissible, then $\cW_{L,T}(\boldsymbol{\gamma}^\star;D,\sigma\mathbf1)=F_{L,T}(\boldsymbol{\gamma}^\star;D,\sigma\mathbf1)=V_{L,T}$ by \eqref{eq:exact}, so that $M_{L,T}\le V_{L,T}$, and \eqref{eq:globalexact} follows from (i). If the worst case of a positive schedule $\boldsymbol{\gamma}'$ is $M_{L,T}$, then $F_{L,T}(\boldsymbol{\gamma}';D,\sigma\mathbf1)\le\cW_{L,T}(\boldsymbol{\gamma}';D,\sigma\mathbf1)=V_{L,T}$, so that $\boldsymbol{\gamma}'$ minimizes the cost, and $\boldsymbol{\gamma}'=\boldsymbol{\gamma}^\star$ by Theorem~\ref{thm:convex}.

(iii) Suppose that $\gamma_0^\star\ge\frac1L$, and let $S_0$ and $S_1$ denote the tails of $\boldsymbol{\gamma}^\star$. By \eqref{eq:gradient}, the condition $\partial_0F=0$ reads
\begin{equation}
 \frac{D^2}{4\sigma^2}\overset{\eqref{eq:gradient},\eqref{eq:kkt}}{=}\frac{\gamma_0^\star S_0^2}{2S_1+\frac1L}\overset{\eqref{eq:tails}}{=}\frac{\gamma_0^\star(\gamma_0^\star+S_1)^2}{2S_1+\frac1L}\ \geq\ \gamma_0^\star\Big(2\gamma_0^\star-\frac1L\Big),\label{eq:firststationarity}
\end{equation}
since the minimum of $\frac{(\gamma_0^\star+S_1)^2}{2S_1+1/L}$ over $S_1\ge0$ is attained at $S_1=\gamma_0^\star-\frac1L$. Solving this quadratic inequality in $L\gamma_0^\star$ gives \eqref{eq:firststepbound}. 
The right-hand side of \eqref{eq:firststepbound} is at most
$\frac32$ if $\sigma\ge\frac{LD}{2\sqrt3}$; in this case,
all the steps are at most $\frac{3}{2L}$, by Theorem~\ref{thm:order}, and the schedule is admissible.
\end{proof}
Theorem~\ref{thm:global} is the composite counterpart of the minimax result of \citet{gs26} for the smooth case, which requires a sufficient noise level as well. A condition that improves with the horizon is given in Corollary~\ref{cor:globalT} below: $\boldsymbol{\gamma}^\star$ is admissible if $T+3\ge\frac{4L^2D^2}{\sigma^2}$. These sufficient conditions are not sharp. Numerically, $\boldsymbol{\gamma}^\star$ is admissible if and only if the ratio $\sigma^2T/(L^2D^2)$ is at least a threshold, which is about $0.239$ for $T=100$ and $0.246$ for $T=400$,
and all its steps are at most $\frac1L$ if and only if this ratio is at least about $0.941$ for $T=100$ and $0.981$ for $T=400$; between these two thresholds, up to about half of its steps are larger than $\frac1L$. When $\boldsymbol{\gamma}^\star$ is not admissible, part (i) locates $M_{L,T}$ between two computable values: the cost of $\boldsymbol{\gamma}^\star$, and the cost of any admissible schedule, for instance the best schedule with steps at most $\frac{3}{2L}$, which is the solution of a convex problem. Proposition~\ref{prop:numericalcertificate}
gives such bounds when the minimizer is known only approximately.

For example, for $T=1$, a step is admissible if and only if $\gamma\le\frac{3}{2L}$, the cost is $\frac{D^2}{4\gamma}+\frac{\sigma^2\gamma}{2(2-L\gamma)}$ if $L\gamma\le1$ and $\frac{D^2}{4\gamma}+\frac{L\sigma^2\gamma^2}{2}$ if $L\gamma\ge1$, and its minimizer and the minimax value are
\begin{align*}
 &\gamma^\star=\frac{2D}{2\sigma+LD}\le\frac1L,&&M_{L,1}=\frac{D\sigma}{2}+\frac{LD^2}{8},&&\text{if }\sigma\ge\tfrac{LD}{2},\\
 &\gamma^\star=\Big(\frac{D^2}{4L\sigma^2}\Big)^{1/3}\in\Big[\frac1L,\frac{3}{2L}\Big],&&M_{L,1}=\frac38\big(4L\sigma^2D^4\big)^{1/3},&&\text{if }\sqrt{\tfrac{2}{27}}\,LD\le\sigma\le\tfrac{LD}{2}.
\end{align*}
If $\sigma<\sqrt{2/27}\,LD$, the minimizer of the cost is not admissible, and part (i), with the step $\frac{3}{2L}$, gives $\frac38(4L\sigma^2D^4)^{1/3}\le M_{L,1}\le\frac{LD^2}{6}+\frac{9\sigma^2}{8L}$.

\subsection{Closed-form solution in the case \texorpdfstring{$\ff=0$}{f=0}}\label{sec:ppa}
In the case $\ff=0$, the cost is the function $F_{0,T}$ of Corollary~\ref{cor:ppa}, which is the exact worst-case value for every positive schedule,
and its minimization can be carried out in closed form. The noise term in \eqref{eq:ppacost} is homogeneous of degree one in the overall scale of the schedule, whereas the initial-error term has degree minus one, so that we can optimize the relative allocations first and the total stepsize second.

\begin{theorem}[The optimal schedule when $\ff=0$]\label{thm:ppaopt}
Define
\begin{equation}
 Q_1\eqdef1,\qquad Q_n\eqdef Q_{n-1}-\tfrac14 Q_{n-1}^2,\qquad n\geq2.\label{eq:qrec}
\end{equation}
Among all positive predetermined $T$-step schedules, the exact minimax value of $\E[\Delta_T]$ in the case $\ff=0$, $V_{0,T}\eqdef\inf_{\boldsymbol{\gamma}\in(0,\infty)^T}F_{0,T}(\boldsymbol{\gamma};D,\sigma\mathbf1)$, is
\begin{equation}
 V_{0,T}=\frac{D\sigma}{2}\sqrt{Q_T}.\label{eq:ppavalue}
\end{equation}
The unique optimal schedule $\widehat{\boldsymbol{\gamma}}$ is generated by
\begin{align}
 \widehat S_0&=\frac{D}{\sigma\sqrt{Q_T}},\label{eq:ppascale}\\
 \widehat\gamma_t&=\frac{Q_{T-t-1}}2\widehat S_t,\qquad
 \widehat S_{t+1}=\widehat S_t-\widehat\gamma_t\quad(0\leq t\leq T-2),\notag\\
 \widehat\gamma_{T-1}&=\widehat S_{T-1}.\label{eq:ppaschedule}
\end{align}
For $T\geq2$, its last two steps are equal, and every earlier step is strictly larger than its successor. For $T=1$, its sole step is $D/\sigma$.
\end{theorem}
\begin{proof}
We first minimize the noise term for a given total stepsize: we prove by induction on $n\geq1$ that, for every $S>0$,
\begin{equation}
 \min_{\substack{\alpha_j>0\ (0\leq j<n)\\\sum_{j=0}^{n-1}\alpha_j=S}}
 \frac1S\sum_{j=0}^{n-1}\frac{\alpha_j^2}{\sum_{k=j}^{n-1}\alpha_k}=Q_n,\label{eq:shape}
\end{equation}
with a unique minimizer. This is clear for $n=1$, since $Q_1=1$. For $n\geq2$, allocate the fraction $\rho\in(0,1)$ of $S$ to the first step. The first term contributes $\rho^2$, and an optimal remaining schedule, of total stepsize $(1-\rho)S$, contributes $(1-\rho)Q_{n-1}$, by the induction hypothesis. Therefore, the minimum in \eqref{eq:shape} is
\begin{align}
 \min_{0<\rho<1}\{\rho^2+(1-\rho)Q_{n-1}\}=Q_{n-1}-\frac{Q_{n-1}^2}{4}\overset{\eqref{eq:qrec}}{=}Q_n.\label{eq:bellman}
\end{align}
The unique minimizer is $\rho=Q_{n-1}/2$, which lies in $(0,1)$ because $0<Q_{n-1}\leq1$. This proves \eqref{eq:shape} and gives the unique relative allocations in \eqref{eq:ppaschedule}.

For a total stepsize $S$, the best value is now
\[
 \min_{\sum\gamma_t=S}F_{0,T}(\boldsymbol{\gamma};D,\sigma\mathbf1)
 \overset{\eqref{eq:ppacost},\eqref{eq:shape}}{=}
 \frac{D^2}{4S}+\frac{\sigma^2SQ_T}{4}.
\]
This strictly convex scalar expression has its unique minimizer at \eqref{eq:ppascale}, with value \eqref{eq:ppavalue}. Exact attainment follows from Corollary~\ref{cor:ppa}.

With two steps remaining, the fraction is $Q_1/2=1/2$, which gives two equal last steps.
With $n\geq3$ steps remaining, the ratio of the next step to the current one is
\begin{equation*}
 \frac{\widehat\gamma_{t+1}}{\widehat\gamma_t}
 \overset{\eqref{eq:ppaschedule}}{=}\frac{Q_{n-2}}{Q_{n-1}}\Big(1-\frac{Q_{n-1}}2\Big)
 \overset{\eqref{eq:qrec}}{=}
 \frac{1-Q_{n-2}/2+Q_{n-2}^2/8}{1-Q_{n-2}/4}<1,
\end{equation*}
because the denominator exceeds the numerator by $Q_{n-2}(2-Q_{n-2})/8>0$.
\end{proof}
Table~\ref{tab:allocations} gives the optimal relative allocations for small horizons. For example, at $T=3$, they are $(3/8,5/16,5/16)$, with $Q_3=39/64$, and the optimal value is $D\sigma\sqrt{39}/16$.
A constant schedule is already strictly suboptimal at this horizon.

\begin{table}[t]\centering\small
\caption{Optimal relative allocations $\widehat\gamma_t/\widehat S_0$ in the case $\ff=0$ for small horizons, compared with the linearly decaying allocations $\gamma_t/S_0$ at the same total stepsize; the equality of the last two optimal steps is visible from $T=2$ on.}\label{tab:allocations}
\begin{tabular}{ccll}\toprule
$T$ & $Q_T$ & Optimal relative allocation & Linear allocation\\\midrule
$1$ & $1$ & $(1)$ & $(1)$\\
$2$ & $3/4$ & $(1/2,1/2)$ & $(2/3,1/3)$\\
$3$ & $39/64$ & $(3/8,5/16,5/16)$ & $(1/2,1/3,1/6)$\\\bottomrule
\end{tabular}
\end{table}

The value $V_{0,T}$ is governed by the sequence $(Q_n)$, for which we have explicit bounds and an asymptotic expansion.
\begin{proposition}[Finite-horizon bounds and the logarithmic correction]\label{prop:qasymptotic}
Let $H_n\eqdef\sum_{j=1}^n1/j$ and $b_n\eqdef4/Q_n$. Then
\begin{align}
 b_1&=4,\qquad b_n\overset{\eqref{eq:qrec}}{=}b_{n-1}+1+\frac1{b_{n-1}-1},\label{eq:brec}\\
 n+3&\overset{\eqref{eq:brec}}{\leq} b_n\overset{\eqref{eq:brec}}{\leq} n+H_{n+1}+\frac32,\label{eq:bbounds}\\
 b_n&\overset{\eqref{eq:brec},\eqref{eq:bbounds}}{=}n+\log n+C_\star+o(1),\qquad C_\star\in\R.\label{eq:basymptotic}
\end{align}
Consequently, $V_{0,T}=D\sigma/\sqrt{T+\log T+C_\star+o(1)}$ and
\begin{equation}
 \frac{D\sigma}{\sqrt{T+H_{T+1}+3/2}}
 \overset{\eqref{eq:ppavalue},\eqref{eq:bbounds}}{\leq} V_{0,T}
 \overset{\eqref{eq:ppavalue},\eqref{eq:bbounds}}{\leq}\frac{D\sigma}{\sqrt{T+3}}.\label{eq:ppabounds}
\end{equation}
\end{proposition}
\begin{proof}
Taking reciprocals in \eqref{eq:qrec} gives \eqref{eq:brec}. Its increments exceed $1$, yielding the lower bound in \eqref{eq:bbounds}.
For $n\geq2$, this gives $b_{n-1}-1\geq n+1$; summing the resulting upper bounds on the increments yields the upper bound in \eqref{eq:bbounds}.
These bounds imply $b_{n-1}-1=n+O(\log n)$, hence
\[
 \frac1{b_{n-1}-1}\overset{\eqref{eq:bbounds}}{=}\frac1n+O\!\left(\frac{\log n}{n^2}\right).
\]
The series of the error terms is absolutely convergent.
Summing \eqref{eq:brec} therefore proves that $b_n-n-\log n$ converges to a finite universal constant $C_\star$, which is \eqref{eq:basymptotic}. Finally, 
$ V_{0,T}\overset{\eqref{eq:ppavalue}}{=}\frac{D\sigma}{\sqrt{b_T}}$,
which gives \eqref{eq:ppabounds} and the claimed expansion.
\end{proof}
The constant $C_\star$ is not new: the sequence $(Q_n)$ is known in optimal stopping theory.
Indeed, $E_n\eqdef1-Q_n/2$ satisfies $E_1=\frac12$ and $E_n\overset{\eqref{eq:qrec}}{=}\frac12\big(1+E_{n-1}^2\big)$, which is the recursion of \citet{mos56} for the largest expected value that can be secured by stopping a sequence of $n$ independent random variables, uniform on $[0,1]$; the optimal relative allocation $\widehat\gamma_t/\widehat S_0$ of Theorem~\ref{thm:ppaopt} is the probability that his optimal rule stops at the draw of index $t$.
Moser deduced the expansion $E_n\sim1-2/(n+\log n+C_\star)$, which is \eqref{eq:basymptotic} \citep{mos56,gm66}, and the value of the constant is $C_\star=1.76799378\ldots$ \citep[Section~5.15]{fin03}.
The convergence in \eqref{eq:basymptotic} is monotone: $b_n-n-\log n$ decreases to $C_\star$, since $b_n-b_{n-1}-1\overset{\eqref{eq:brec}}{=}\frac1{b_{n-1}-1}\overset{\eqref{eq:bbounds}}{\le}\frac1{n+1}<\log\frac{n}{n-1}$ for every $n\ge2$, and the proof of Proposition~\ref{prop:qasymptotic} shows that the difference is $O(\log n/n)$.
The bounds \eqref{eq:ppabounds} locate $V_{0,T}$ for every finite horizon, and the expansion shows that the optimal schedule behaves as if the horizon were extended by $\log T+C_\star$ iterations, in the sense that $V_{0,T}$ is the leading-order rate $D\sigma/\sqrt T$ evaluated at this effective horizon.
For a constant stepsize in the case $\ff=0$, the cost is
\begin{equation}
 F_{0,T}(\gamma\mathbf1;D,\sigma\mathbf1)
 \overset{\eqref{eq:ppacost},\eqref{eq:tails}}{=}
 \frac{D^2}{4T\gamma}+\frac{\sigma^2\gamma H_T}{4}.
 \label{eq:constantppacost}
\end{equation}
Optimizing this scalar expression gives
\begin{equation*}
 \min_{\gamma>0}F_{0,T}(\gamma\mathbf1;D,\sigma\mathbf1)
 \overset{\eqref{eq:constantppacost}}{=}\frac{D\sigma}{2}\sqrt{\frac{H_T}{T}},\qquad
 \gamma=\frac{D}{\sigma\sqrt{TH_T}}.
\end{equation*}
Thus, optimizing the whole schedule, instead of using the best constant stepsize, removes a factor
asymptotic to $\frac12\sqrt{\log T}$ from the worst-case value, without changing the algorithm.

\subsection{Nonuniform noise: exact design and when stopping is optimal}\label{sec:nonuniform}
In this subsection, the variance bound is a prescribed profile $\boldsymbol{\sigma}=(\sigma_t)_{t=0}^{T-1}$, no longer uniform. The uniformity was essential to the positivity of the optimizer in Theorem~\ref{thm:convex}; in the case $\ff=0$, the nonuniform problem can be solved completely, and its solution may stop before the horizon. Here, and only here, we allow zero steps, $\gamma_t=0$, interpreted as $x^{t+1}=x^t$, and we require a positive total stepsize $S_0>0$; the summands of \eqref{eq:ppacost} with a zero tail $S_t=0$ are set to zero. Since a zero step can be removed, with the horizon reduced accordingly, this extended cost is still the exact worst-case value, by Corollary~\ref{cor:ppa}.

\begin{theorem}[Design for a positive variance profile]\label{thm:heterogeneous}
Let $D>0$ and $\sigma_t>0$ for all $t$. Define $A_{T-1}\eqdef\sigma_{T-1}^2$, $\widehat\rho_{T-1}\eqdef1$, and, backward for $t=T-2,\ldots,0$,
\begin{align}
 \widehat\rho_t&\eqdef\min\left\{\frac{A_{t+1}}{2\sigma_t^2},1\right\},\label{eq:heteroratio}\\
 A_t&\eqdef\sigma_t^2\widehat\rho_t^{\,2}+(1-\widehat\rho_t)A_{t+1}
 \overset{\eqref{eq:heteroratio}}{=}\begin{cases}
 A_{t+1}-A_{t+1}^2/(4\sigma_t^2),&A_{t+1}\leq2\sigma_t^2,\\
 \sigma_t^2,&A_{t+1}\geq2\sigma_t^2.
 \end{cases}\label{eq:heterorec}
\end{align}
Then $A_t>0$. The unique optimal nonnegative schedule is
\begin{equation}
 \widehat S_0=\frac D{\sqrt{A_0}},\qquad
 \widehat\gamma_t=\widehat\rho_t \widehat S_t,\qquad
 \widehat S_{t+1}=(1-\widehat\rho_t)\widehat S_t,
 \label{eq:heteroschedule}
\end{equation}
and its value is $D\sqrt{A_0}/2$. Its positive steps form a prefix. The infimum over strictly positive schedules is the same, and is attained if and only if
\begin{equation}
 A_{t+1}<2\sigma_t^2\qquad(0\leq t<T-1).\label{eq:positiveprofile}
\end{equation}
Both the backward recursion and the forward schedule construction use $O(T)$ arithmetic operations.
\end{theorem}
\begin{proof}
For a remaining total $S$, allocate $\rho S$ to its first step. If the optimal remaining noise cost per unit total is $A_{t+1}$, the current cost per unit total is $\sigma_t^2\rho^2+(1-\rho)A_{t+1}$. This strictly convex quadratic has the unique minimizer \eqref{eq:heteroratio} on $[0,1]$, with value \eqref{eq:heterorec}. The last stage has cost $\sigma_{T-1}^2S$, starting the induction. Positivity follows since $\sigma_t^2>0$ and $A_{t+1}>0$.

For a direct global certificate, define
\begin{equation}
 \delta_t(\rho)\eqdef \sigma_t^2(\rho-\widehat\rho_t)^2
   +\mathbf1_{\{\widehat\rho_t=1\}}(A_{t+1}-2\sigma_t^2)(1-\rho),
 \qquad t<T-1.\label{eq:heterogap}
\end{equation}
Completing the square gives, for $0\leq\rho\leq1$,
\begin{equation}
 \sigma_t^2\rho^2+(1-\rho)A_{t+1}-A_t
 \overset{\eqref{eq:heteroratio},\eqref{eq:heterorec},\eqref{eq:heterogap}}{=}
 \delta_t(\rho)\overset{\eqref{eq:heterogap},\eqref{eq:heteroratio}}{\geq}0.
 \label{eq:bellmangap}
\end{equation}
For any nonnegative schedule with $S_0>0$, set $\rho_t\eqdef\gamma_t/S_t$ when $S_t>0$; interpret each term with $S_t=0$ as zero. Multiplying by $S_t$ and telescoping proves the exact identity
\begin{align}
 F_{0,T}(\boldsymbol{\gamma};D,\boldsymbol{\sigma})-\frac{D\sqrt{A_0}}2
 &\overset{\eqref{eq:ppacost},\eqref{eq:tails},\eqref{eq:bellmangap}}{=}
 \frac{(D-S_0\sqrt{A_0})^2}{4S_0}
 +\frac14\sum_{t=0}^{T-2}S_t\delta_t(\rho_t).
 \label{eq:profilecertificate}
\end{align}
Every term is nonnegative. Equality forces the total in \eqref{eq:heteroschedule} and each active ratio to equal $\widehat\rho_t$, proving optimality and uniqueness even when a tail is zero. Since every ratio is positive, there are no interior zero steps before stopping. A ratio equals one precisely at the threshold excluded by \eqref{eq:positiveprofile}. Finally, any nonnegative schedule is the limit of positive ones with the same total, and the extended cost is continuous on that simplex: $0\leq\gamma_t^2/S_t\leq\gamma_t$. This proves the infimum and attainment assertions.
\end{proof}
For example, with $T=2$ and $\boldsymbol{\sigma}^2=(1,4)$,
the recursion gives $A_1=4$, $\widehat\rho_0=\min\big\{\frac42,1\big\}=1$, and $A_0=1$,
so that the optimal schedule is $(D,0)$, with value $\frac D2$: the second step, four times noisier than the first, would cost more than the reduction of the tail it buys, and it is better to stop after one step. Requiring both steps to be positive preserves the infimum but loses attainment. The variance profile is known in advance, so that no adaptivity is involved.

If some $\sigma_j=0$ and $D>0$, the infimum over positive schedules is zero and is not attained. To see this, make the noiseless step arbitrarily large and all other steps arbitrarily small in \eqref{eq:ppacost}. The initial-error term vanishes, earlier noise terms have a large denominator, and later ones are bounded by $\sigma_t^2\gamma_t/4$. Every finite positive schedule still has a positive initial-error term. This separates the zero-variance case from the positive-profile theorem.

\begin{corollary}[An exact certificate for the schedule error when $\ff=0$]\label{cor:schedulegap}
Let $\ff=0$ and $D,\sigma>0$, and let $\boldsymbol{\gamma}$ be any positive schedule; put
$\rho_t\eqdef\gamma_t/S_t$. Then
\begin{align}
 F_{0,T}(\boldsymbol{\gamma};D,\sigma\mathbf1)-V_{0,T}
 &\overset{\eqref{eq:profilecertificate},\eqref{eq:heterorec},\eqref{eq:qrec},\eqref{eq:ppavalue}}{=}
 \frac{(D-\sigma S_0\sqrt{Q_T})^2}{4S_0}
 +\frac{\sigma^2}{4}\sum_{t=0}^{T-2}S_t
 \left(\rho_t-\frac{Q_{T-t-1}}2\right)^2.
 \label{eq:uniformgap}
\end{align}
In particular, multiplying the optimal schedule by $\alpha>0$ changes its value by the exact factor
\begin{equation}
 \frac{F_{0,T}(\alpha\widehat{\boldsymbol{\gamma}};D,\sigma\mathbf1)}{V_{0,T}}
 \overset{\eqref{eq:uniformgap},\eqref{eq:ppascale},\eqref{eq:ppavalue}}{=}
 \frac{\alpha+\alpha^{-1}}2.
 \label{eq:scalesensitivity}
\end{equation}
\end{corollary}
\begin{proof}
For a uniform variance bound, the backward recursion gives
$A_t=\sigma^2Q_{T-t}$ and $\widehat\rho_t=Q_{T-t-1}/2<1$ before the last stage. Substitution into \eqref{eq:profilecertificate} proves \eqref{eq:uniformgap}. Scaling preserves every ratio; substituting $S_0=\alpha D/(\sigma\sqrt{Q_T})$ gives \eqref{eq:scalesensitivity}.
\end{proof}
The first term of \eqref{eq:uniformgap} is a scale error and the second an allocation error, so that the optimality of $\widehat{\boldsymbol{\gamma}}$ has an explicit quantitative content. The scale factor \eqref{eq:scalesensitivity} also describes the effect of a misspecified ratio $D/\sigma$ when the horizon-dependent shape is correct: the schedule is not parameter-free.

\section{Asymptotic minimax law}\label{sec:asymptotic}

The minimax value $M_{L,T}$ of \eqref{eq:globalvalue} is known exactly when the minimizer of the cost is admissible (Theorem~\ref{thm:global}). This section determines its behavior as $T\to\infty$, for every fixed noise level $\sigma>0$, and shows that this admissibility holds when the horizon is long enough. Recall that $\dim\cH\ge T+1$ is assumed for every horizon $T$ considered; in particular, $\cH$ is infinite-dimensional in the statements where $T\to\infty$.
The denominators \eqref{eq:denominators} satisfy $d_t\le2S_t$, so that $F_{L,T}\ge F_{0,T}$ at every positive schedule; taking the infimum over all positive schedules and using Theorem~\ref{thm:global} gives, for every $T\ge1$,
\begin{equation}
 V_{0,T}\overset{\eqref{eq:ppacost},\eqref{eq:cost},\eqref{eq:optvalue}}{\leq}
 V_{L,T}\overset{\eqref{eq:universallower},\eqref{eq:optvalue},\eqref{eq:globalvalue}}{\leq}
 M_{L,T}.
 \label{eq:globalorder}
\end{equation}

\subsection{The optimal schedule of the case \texorpdfstring{$\ff=0$}{f=0} transfers the full logarithmic correction}\label{sec:transfer}
Recall that $\widehat{\boldsymbol{\gamma}}$ denotes the optimal schedule of the case $\ff=0$ (Theorem~\ref{thm:ppaopt}). Its shape controls every local noise contribution tightly enough that adding curvature changes its value only at order $T^{-2}$, so that it carries the second-order law of Proposition~\ref{prop:qasymptotic} over to every $L>0$.

\begin{theorem}[A universal second-order minimax law]\label{thm:secondorder}
Let $L,D,\sigma>0$, and suppose that $LD\sqrt{Q_T}\leq\sigma$. Define
\begin{equation}
 \varepsilon_T\eqdef \frac{LDQ_T^{3/2}}{2\sigma}\leq\frac12.
 \label{eq:transferparameter}
\end{equation}
The schedule $\widehat{\boldsymbol{\gamma}}$ satisfies $\widehat\gamma_t\leq\frac1L$. In addition to \eqref{eq:globalorder},
\begin{align}
 M_{L,T}
 &\overset{\eqref{eq:globalvalue}}{\leq}\cW_{L,T}(\widehat{\boldsymbol{\gamma}};D,\sigma\mathbf1)\leq V_{0,T}+\frac{LD^2Q_T^2}{8(1-\varepsilon_T)}.
 \label{eq:secondordersandwich}
\end{align}
For $Z_T\in\{V_{L,T},M_{L,T},\cW_{L,T}(\widehat{\boldsymbol{\gamma}};D,\sigma\mathbf1)\}$, these bounds imply
\begin{equation}
 0\overset{\eqref{eq:globalorder},\eqref{eq:secondordersandwich},\eqref{eq:ppavalue}}{\leq} b_T-\left(\frac{D\sigma}{Z_T}\right)^2
 \overset{\eqref{eq:globalorder},\eqref{eq:secondordersandwich},\eqref{eq:ppavalue}}{\leq}\frac{2LD\sqrt{Q_T}}{\sigma(1-\varepsilon_T)}.
 \label{eq:effectivehorizon}
\end{equation}
Consequently, for each fixed $L,D,\sigma>0$,
\begin{equation}
 \left(\frac{D\sigma}{Z_T}\right)^2-T-\log T
 \overset{\eqref{eq:effectivehorizon},\eqref{eq:basymptotic}}{\longrightarrow} C_\star,
 \label{eq:universalconstant}
\end{equation}
with the same universal constant as in \eqref{eq:basymptotic}. Equivalently,
$Z_T=D\sigma/\sqrt{T+\log T+C_\star+o(1)}$.
\end{theorem}
The proof is given in Appendix~\ref{app:proofsecondorder}. The lower comparison is \eqref{eq:globalorder}, and the upper bound is obtained by evaluating the cost at $\widehat{\boldsymbol{\gamma}}$. Indeed, for every schedule with steps at most $\frac1L$, \eqref{eq:convexcomposition} gives the exact curvature penalty
\begin{equation}
 F_{L,T}-F_{0,T}\overset{\eqref{eq:cost},\eqref{eq:ppacost},\eqref{eq:convexcomposition}}{=}
 \frac{\sigma^2L}{2}\sum_t\frac{r_t^2}{1-Lr_t},
 \label{eq:generalpenalty}
\end{equation}
and the numbers $\widehat r_t=\widehat\gamma_t^2/(2\widehat S_t)$ of the schedule $\widehat{\boldsymbol{\gamma}}$ are small: their sum is $\frac{D\sqrt{Q_T}}{2\sigma}$ and each of them is at most $\frac{DQ_T^{3/2}}{2\sigma}$.

The new information in Theorem~\ref{thm:secondorder} is below the leading constant:
curvature can alter neither the logarithmic correction nor its limiting constant,
even when the competing schedules may take arbitrarily large predetermined steps.
The comparison is between schedules of the prescribed algorithm, and the limit is taken for fixed $L,D,\sigma$; it is not uniform when these parameters vary with $T$.

This is where the composite class differs from the smooth one at second order. In the smooth case, the minimax rate of \citet{gs26} recalled in Section~\ref{sec:classical}, $r_\star=\frac{\sigma D}{\sqrt T}-\frac{\sigma^{3/2}\sqrt D}{\sqrt L\,T^{3/4}}+o(T^{-3/4})$ under their sufficient-variance condition, has a second-order term that depends on $L$ and improves on the leading term at order $T^{-3/4}$; written as an effective horizon, $(\sigma D/r_\star)^2=T+2\sqrt{\sigma/(LD)}\,T^{3/4}+o(T^{3/4})$. In the composite case, the effective horizon is $T+\log T+C_\star+o(1)$,
whatever $L$:
the prox removes the gain from curvature at second order, just as it removed the $2/L$ from the initial-error term, and the two classes, which share the leading term $\frac{\sigma D}{\sqrt T}$, separate at order $T^{-3/4}$.

The bounds of Theorem~\ref{thm:secondorder} also control the largest step of the minimizer of the cost, and hence its admissibility, by a condition that improves with the horizon.
\begin{corollary}[Minimax optimal schedules for long horizons]\label{cor:globalT}
Let $L,D,\sigma>0$ and let $\boldsymbol{\gamma}^\star$ be the minimizer of the cost over all positive schedules, as in Theorem~\ref{thm:global}. If $LD\sqrt{Q_T}\le\sigma$, which is the case if $T+3\ge\frac{4L^2D^2}{\sigma^2}$, then all the steps of $\boldsymbol{\gamma}^\star$ are smaller than $\frac{1}{L}$; hence $\boldsymbol{\gamma}^\star$ is admissible, $M_{L,T}=V_{L,T}$, and $\boldsymbol{\gamma}^\star$ is the unique minimax optimal schedule.
\end{corollary}
\begin{proof}
Let $S_i$ and $d_i$ denote the tails and the denominators \eqref{eq:denominators} of $\boldsymbol{\gamma}^\star$. Summing the conditions \eqref{eq:kkt}, multiplied by $\gamma_k^\star$, and using \eqref{eq:gradient} gives
\begin{align*}
 \frac{D^2}{4S_0}&\overset{\eqref{eq:gradient},\eqref{eq:kkt}}{=}\sigma^2\sum_{i=0}^{T-1}\frac{(\gamma_i^\star)^2\big(S_{i+1}+\min\{\gamma_i^\star,\frac1L\}\big)}{d_i^2}\\
 &\overset{\eqref{eq:denominators},\eqref{eq:omega}}{\leq}\sigma^2\sum_{i=0}^{T-1}\frac{(\gamma_i^\star)^2}{d_i}\overset{\eqref{eq:cost},\eqref{eq:optvalue}}{=}2\Big(V_{L,T}-\frac{D^2}{4S_0}\Big),
\end{align*}
where the inequality holds because $\min\{\gamma,\frac1L\}\le\omega(\gamma)$ for every $\gamma>0$, so that $S_{i+1}+\min\{\gamma_i^\star,\frac1L\}\le d_i$. Hence, the initial-error term $\frac{D^2}{4S_0}$ is at most $\frac23V_{L,T}$.
Suppose that $\gamma_0^\star\ge\frac1L$. By \eqref{eq:firststationarity}, $\frac{D^2}{4S_0}=\frac{\sigma^2\gamma_0^\star(\gamma_0^\star+S_1)}{2S_1+1/L}>\frac{\sigma^2\gamma_0^\star}{2}$, because $2\gamma_0^\star>\frac1L$, so that $\gamma_0^\star<\frac{4V_{L,T}}{3\sigma^2}$. Since $LD\sqrt{Q_T}\le\sigma$, Theorem~\ref{thm:secondorder} applies, with $\varepsilon_T\le\frac12$, and gives
\begin{equation*}
 V_{L,T}\overset{\eqref{eq:globalorder},\eqref{eq:secondordersandwich},\eqref{eq:ppavalue}}{\leq}\frac{D\sigma\sqrt{Q_T}}{2}+\frac{LD^2Q_T^2}{4}
 \leq\frac{D\sigma\sqrt{Q_T}}{2}+\frac{D\sigma Q_T^{3/2}}{4}\leq\frac34D\sigma\sqrt{Q_T},
\end{equation*}
so that $L\gamma_0^\star<\frac{LD\sqrt{Q_T}}{\sigma}\le1$, a contradiction. Hence $\gamma_0^\star<\frac1L$, all the steps of $\boldsymbol{\gamma}^\star$ are smaller than $\frac1L$, by Theorem~\ref{thm:order}, and Theorem~\ref{thm:global} applies.
Finally, $Q_T=4/b_T\le4/(T+3)$ by \eqref{eq:bbounds}.
\end{proof}
For fixed $L,D,\sigma$, the minimizer of the cost, computed by the scalar construction of Theorem~\ref{thm:scalarconstruction}, 
is thus \emph{the} minimax optimal schedule for every $T\ge\frac{4L^2D^2}{\sigma^2}-3$.

\subsection{Linear decay: a simple leading-order alternative}
\begin{corollary}[An explicit asymptotically minimax schedule]\label{thm:asymptotic}
Fix $L,D,\sigma>0$. For all $T$ such that $LD\leq\sigma\sqrt{T+1}$, set
\begin{equation}
 a\eqdef\frac{D}{\sigma T\sqrt{T+1}},\qquad \kappa\eqdef\frac{LaT}{T+1},\qquad
 \gamma_t\eqdef a(T-t).\label{eq:linearschedule}
\end{equation}
Together with \eqref{eq:globalorder}, the following upper bound holds:
\begin{align}
 M_{L,T}
 &\leq\frac{D\sigma}{\sqrt{T+1}}
 \left(1-\frac{H_{T+1}-1}{2T}\right)
 +\frac{LD^2}{2T(T+1)(1-\kappa)}\label{eq:sandwich}\\
 &\leq\frac{D\sigma}{\sqrt{T+1}}
 +\frac{LD^2}{2T(T+1)(1-\kappa)}.\notag
\end{align}
Here $\kappa\leq1/(T+1)<1$. For fixed $L,D,\sigma$,
\begin{equation}
 \lim_{T\to\infty}\frac{\sqrt T\,M_{L,T}}{D\sigma}
 \overset{\eqref{eq:universalconstant}}{=}
 \lim_{T\to\infty}\frac{\sqrt T\,V_{L,T}}{D\sigma}
 \overset{\eqref{eq:universalconstant}}{=}1.\label{eq:minimaxlimit}
\end{equation}
The schedule \eqref{eq:linearschedule} itself attains the same limiting constant.
\end{corollary}
\begin{proof}
The lower comparisons are those of \eqref{eq:globalorder}. The hypothesis ensures that $\gamma_0=D/(\sigma\sqrt{T+1})\leq\frac1L$, so that \eqref{eq:exact} applies to the proposed schedule.
Write $n\eqdef T-t$. Then
\begin{equation}
 S_t\overset{\eqref{eq:tails},\eqref{eq:linearschedule}}{=}\frac{an(n+1)}2.\label{eq:lineartail}
\end{equation}
Hence
\begin{equation}
 r_t\eqdef \frac{\gamma_t^2}{2S_t}
 \overset{\eqref{eq:linearschedule},\eqref{eq:lineartail}}{=}\frac{an}{n+1}\leq\frac{aT}{T+1}\leq a.\label{eq:linearsums}
\end{equation}
Since $Lr_t\leq\kappa<1$, the curvature penalty satisfies
\begin{align}
 F_{L,T}-F_{0,T}
 &\overset{\eqref{eq:cost},\eqref{eq:ppacost},\eqref{eq:convexcomposition}}{=}
 \frac{\sigma^2L}{2}\sum_{t=0}^{T-1}\frac{r_t^2}{1-Lr_t}\overset{\eqref{eq:linearsums},\eqref{eq:linearschedule}}{\leq}
 \frac{LD^2}{2T(T+1)(1-\kappa)}.\label{eq:curvaturepenalty}
\end{align}
The zero-curvature part can be evaluated exactly, using $\sum_{n=1}^T n/(n+1)=T+1-H_{T+1}$:
\begin{equation}
 F_{0,T}(\boldsymbol{\gamma};D,\sigma\mathbf1)
 \overset{\eqref{eq:ppacost},\eqref{eq:linearschedule},\eqref{eq:linearsums}}{=}
 \frac{D\sigma}{\sqrt{T+1}}\left(1-\frac{H_{T+1}-1}{2T}\right).\label{eq:linearcost}
\end{equation}
Combining \eqref{eq:curvaturepenalty} and \eqref{eq:linearcost} proves \eqref{eq:sandwich}. For fixed positive parameters, $\kappa=O(T^{-3/2})$, and the penalty in \eqref{eq:curvaturepenalty} is $O(LD^2/T^2)$. Proposition~\ref{prop:qasymptotic} together with \eqref{eq:globalorder} and \eqref{eq:sandwich} proves \eqref{eq:minimaxlimit} and the claim for the explicit schedule.
\end{proof}
The order $T^{-1/2}$ with linearly decaying stepsizes was known \citep[Theorems 3.4 and 3.5]{lz24},
and the leading constant $1$ is that of the smooth minimax rate \citep{gs26}.
What Corollary~\ref{thm:asymptotic} adds is that this constant is attainable and asymptotically unavoidable for the whole composite class, 
even against predetermined schedules with arbitrarily large steps,
and that an explicit, simple schedule achieves it. This asymptotic statement alone does not settle the exact finite-horizon value $M_{L,T}$,
and nothing is said about adaptive, accelerated, or averaged methods.

\begin{figure}[t]\centering
\includegraphics[width=.72\linewidth]{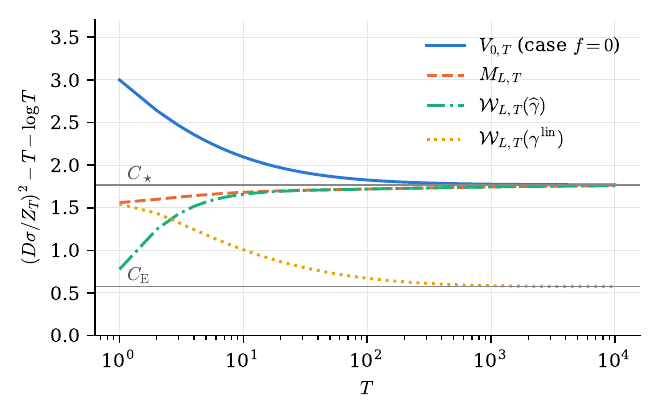}
\caption{Effective horizon $(D\sigma/Z_T)^2-T-\log T$ for $L=D=\sigma=1$, as a function of $T$, for the exact minimax value $V_{0,T}$ of the case $\ff=0$ (Theorem~\ref{thm:ppaopt}), 
the minimax value $M_{L,T}=V_{L,T}$ (Theorems~\ref{thm:scalarconstruction} and~\ref{thm:global}),
the exact worst case of the schedule $\widehat{\boldsymbol{\gamma}}$ of the case $\ff=0$ (Theorem~\ref{thm:secondorder}), and the exact worst case of linear decay (Corollary~\ref{thm:asymptotic}). 
The first three curves converge to $C_\star\approx1.768$; 
the second one and the third one cannot be distinguished for $T\ge10$.
Linear decay converges to Euler's constant $C_{\rm E}\approx0.577$.
}\label{fig:horizon}
\end{figure}

The distinction between leading-order and finer optimality is quantitative. For $T\ge2$, the last two steps of \eqref{eq:linearschedule} are $2a$ and $a$, with tail $S_{T-2}=3a$ and first fraction $\frac23$, whereas the optimal fraction with two steps remaining is $\frac{Q_1}2=\frac12$ in the case $\ff=0$. Keeping only this term in the identity \eqref{eq:uniformgap} gives
\begin{equation}
 F_{0,T}(\boldsymbol{\gamma};D,\sigma\mathbf1)-V_{0,T}
 \overset{\eqref{eq:uniformgap},\eqref{eq:linearschedule}}{\geq}
 \frac{\sigma^2}{4}(3a)\left(\frac23-\frac12\right)^2
 \overset{\eqref{eq:linearschedule}}{=}\frac{D\sigma}{48T\sqrt{T+1}}.
 \label{eq:lineardefect}
\end{equation}
Thus the linear schedule loses at least a term of order $T^{-3/2}$ in the case $\ff=0$, whereas the optimal schedule of that case loses at most a term of order $T^{-2}$ when transferred to the composite class.
More precisely, for fixed $L,D,\sigma>0$,
\begin{equation}
 \liminf_{T\to\infty}\frac{T^{3/2}}{D\sigma}
 [\cW_{L,T}(\boldsymbol{\gamma}^{\rm lin};D,\sigma\mathbf1)-M_{L,T}]
 \overset{\eqref{eq:exact},\eqref{eq:generalpenalty},\eqref{eq:secondordersandwich},\eqref{eq:lineardefect}}{\geq}\frac1{48},
 \label{eq:linearseparation}
\end{equation}
where $\boldsymbol{\gamma}^{\rm lin}$ denotes the schedule \eqref{eq:linearschedule}. The curvature-transfer error is $O(T^{-2})$, so that it vanishes after this normalization.
We can identify exactly which constant linear decay misses. 
Let $C_{\rm E}\eqdef\lim_{n\to\infty}(H_n-\log n)$ be Euler's constant.
Expanding $(1-z)^{-2}=1+2z+O(z^2)$ at $z=(H_{T+1}-1)/(2T)$ gives
\begin{align}
 \left(\frac{D\sigma}{F_{0,T}(\boldsymbol{\gamma}^{\rm lin})}\right)^2
 &\overset{\eqref{eq:linearcost}}{=}
 \frac{T+1}{[1-(H_{T+1}-1)/(2T)]^2}\notag\\
 &=T+H_{T+1}+O((\log T)^2/T).
 \label{eq:lineareffective}
\end{align}
By \eqref{eq:curvaturepenalty}, replacing $F_{0,T}(\boldsymbol{\gamma}^{\rm lin})$ by the exact composite value changes this expression by $O(T^{-1/2})$. Therefore
\begin{equation}
 \left(\frac{D\sigma}{\cW_{L,T}(\boldsymbol{\gamma}^{\rm lin};D,\sigma\mathbf1)}\right)^2
 -T-\log T
 \overset{\eqref{eq:exact},\eqref{eq:curvaturepenalty},\eqref{eq:lineareffective}}{\longrightarrow}
 C_{\rm E}.
 \label{eq:linearconstant}
\end{equation}
Expanding the reciprocal square roots in \eqref{eq:universalconstant} and \eqref{eq:linearconstant} 
shows that the normalized gap in \eqref{eq:linearseparation} has limit $(C_\star-C_{\rm E})/2$. Hence $C_\star\geq C_{\rm E}+1/24$.
In fact, $C_\star-C_{\rm E}\approx1.191$.
Linear decay and the optimal schedule agree in their leading constant and logarithmic term, but differ in the next constant.
Figure~\ref{fig:horizon} illustrates these asymptotics for $L=D=\sigma=1$, by plotting the effective horizon $(D\sigma/Z_T)^2-T-\log T$ against $T$ for the four values 
$Z_T\in\{V_{0,T},M_{L,T},\cW_{L,T}(\widehat{\boldsymbol{\gamma}};D,\sigma\mathbf1),\cW_{L,T}(\boldsymbol{\gamma}^{\rm lin};D,\sigma\mathbf1)\}$, all of which are exact worst-case values here.
The first three converge to $C_\star$: the first from above, at the rate $O(\log T/T)$, as noted in Section~\ref{sec:ppa}; the other two from below, as the figure shows, within the $O(T^{-1/2})$ margin of \eqref{eq:effectivehorizon}.
The costs of the optimal schedule and of $\widehat{\boldsymbol{\gamma}}$ differ by $20\%$ at $T=1$ and by $4.4\%$ at $T=2$, but by less than $10^{-3}$ in relative value for $T\ge10$, and by less than $10^{-6}$ for $T\ge100$, where the two curves cannot be distinguished. 
Here $\sigma\ge\frac{LD}{2\sqrt3}$, so that $M_{L,T}=V_{L,T}$ for every $T$, by Theorem~\ref{thm:global}:
the second curve is the exact minimax value over all positive predetermined schedules.
Linear decay converges to $C_{\rm E}$, as \eqref{eq:linearconstant} states.

\section{Conclusion}\label{sec:conclusion}
We have obtained the exact worst-case value of the expected last-iterate objective gap of \algno for every admissible predetermined schedule; this class contains all the schedules with steps at most $\frac{3}{2L}$ and allows early steps up to nearly $\frac2L$. The formula separates an initial-error term, $\frac{D^2}{4S_0}$, from one contribution per noise injection, $\frac{\sigma_t^2\gamma_t^2}{2d_t}$, in which the remaining stepsize sum is the resource that damps the injection made at iteration $t$; the prox makes this contribution larger than in the smooth case
for the steps smaller than $\frac1L$, and only for them.
This structure explains the design results of Sections~\ref{sec:design} and~\ref{sec:asymptotic}: optimal schedules decrease, since an early injection is damped by a long tail and a late one is not; under a nonuniform variance profile, a sufficiently noisy future is best discarded by stopping; and matching the whole schedule, not only the order of its decay, matters at second order,
where linear decay and the optimal schedule part ways.
Moreover, the minimizer of the cost over all positive schedules is exactly
minimax optimal among all positive predetermined schedules whenever it is admissible, which is the case unless the noise level is low: $\sigma\ge\frac{LD}{2\sqrt3}$ or $T+3\ge\frac{4L^2D^2}{\sigma^2}$ suffices.
In the case $\ff=0$, the exact value is known for every schedule, the optimal schedule is explicit, and it transfers to every $L>0$ at the cost of an $O(T^{-2})$ curvature penalty, which yields the universal law $(D\sigma/M_{L,T})^2=T+\log T+C_\star+o(1)$ and the contrast with the smooth class, which gains from curvature at order $T^{-3/4}$.

Several questions are open and left for future work. First, the schedules that are not admissible: the cost $F_{L,T}$ is then a lower bound only, and the exact deterministic rates recalled in Section~\ref{sec:classical} show that other regimes, governed by quadratic instances, appear for larger stepsizes. 
Second, the finite-horizon minimax problem at low noise levels, when the minimizer of the cost is not admissible, which happens, numerically, for long horizons, when $\sigma^2T$ is smaller than about $\frac{L^2D^2}{4}$: the minimax value is then only located between two computable values (Theorem~\ref{thm:global}).
The acceleration of \algn{Prox-GD} by large stepsizes \citep{ba25} belongs to the deterministic limit of this question. Third, adaptive schedules, which break the predictability on which the analysis relies, and other noise models \citep{kha23,vfa26}, 
which call for other arguments.

\appendix
\section{Complete interpolation certificate}\label{app:certificate}
This appendix proves Lemma~\ref{lem:local} and records enough algebra to reproduce every coefficient.
Let $\gamma$, $\tau$, $\theta$, $\tau'$, $c$, and $d$ be as in Lemma~\ref{lem:local}, and let $x,e\in\cH$. In addition to the output $q$ of the noisy step, defined in \eqref{eq:local}, we use the points
\begin{equation}
 p\eqdef\prox_{\tau\Psi}(x),\qquad
 \qc\eqdef\prox_{\gamma \RR}\big(x-\gamma\nabla \ff(x)\big),\label{eq:localpoints}
\end{equation}
where $\qc$ is the output of the noiseless step. With $\gc\eqdef\min\big\{\gamma,\frac1L\big\}$, so that $\omega(\gamma)=\gc(2-L\gc)$, and $\mu\eqdef\frac{L(\gamma-\gc)}{1+L\tau'}$, the vector $B(x)$ of Lemma~\ref{lem:local} is
\begin{equation}
 B(x)\eqdef-\gamma\frac{x-p}{\tau}
 -\frac{\gamma\gc(\tau-\gamma)}{\tau d}
       \big(\nabla \ff(p)-\nabla \ff(x)\big)
 -\mu\Big(\gamma\frac{x-p}{\tau}-(x-\qc)\Big).\label{eq:Bc}
\end{equation}
It does not depend on $e$.
When $L\gamma\le1$, we have $\mu=0$, and the point $\qc$ is not involved.
The mapping $B$ is Lipschitz continuous, because $\nabla\ff$ is Lipschitz continuous and the proxes in \eqref{eq:localpoints} are nonexpansive.
The case $L\gamma\le1$ is treated first; the modifications for $L\gamma>1$ are given in Section~\ref{app:long}.
For every $y,z\in\cH$ and $s_z\in\partial \RR(z)$, define the interpolation residuals
\begin{align}
 I_f(y,z)&\eqdef\ff(y)-\ff(z)-\ip{\nabla \ff(z)}{y-z}
              -\tfrac1{2L}\norm{\nabla \ff(y)-\nabla \ff(z)}^2\overset{\eqref{eq:interp}}{\geq}0,\label{eq:If}\\
 I_R(y,z;s_z)&\eqdef\RR(y)-\RR(z)-\ip{s_z}{y-z}\overset{\eqref{eq:subgradient}}{\geq}0.\label{eq:Ig}
\end{align}
For completeness, we recall that the first inequality follows by applying the descent lemma to
$\ell(w)=\ff(w)-\ip{\nabla \ff(z)}w$.
By convexity, $z$ is a minimizer of $\ell$, and one gradient step from $y$ gives
$\ell(z)\leq \ell(y)-\norm{\nabla \ell(y)}^2/(2L)$, which is \eqref{eq:If}.
The second inequality is the subgradient inequality.

For $\tau'>0$, put $u\eqdef\prox_{\tau'\Psi}(q)$; at $\tau'=0$, set $u\eqdef q$. Write $g_x\eqdef\nabla \ff(x)$, $g_p\eqdef\nabla \ff(p)$, and $g_u\eqdef\nabla \ff(u)$. Select
\begin{align}
 s_p&\eqdef(x-p)/\tau-g_p\overset{\eqref{eq:localpoints}}{\in}\partial \RR(p),\qquad
 s_q\eqdef(x-q)/\gamma-g_x-e\overset{\eqref{eq:local}}{\in}\partial \RR(q),\notag\\
 s_u&\eqdef(q-u)/\tau'-g_u\in\partial \RR(u)\quad(\tau'>0),
 \qquad s_u\eqdef s_q\quad(\tau'=0).\label{eq:subgradients}
\end{align}
These inclusions follow from proximal optimality; $\partial\Psi=\nabla \ff+\partial \RR$ because $\ff$ is finite and continuous. Introduce
\begin{align}
 G_p&\eqdef g_p+s_p,\quad G_x\eqdef g_x+e+s_q,\quad U\eqdef g_u+s_u,\notag\\
 a&\eqdef G_p-U,\quad b\eqdef G_x-U,\quad A_f\eqdef g_p-g_u,\quad B_f\eqdef g_x-g_u.\label{eq:vectors}
\end{align}
The mixed vector $G_x$ is not asserted to be a subgradient at $x$. The exact displacement relations are
\begin{align}
 x-p&\overset{\eqref{eq:subgradients},\eqref{eq:vectors}}{=}\tau G_p,
 &x-q&\overset{\eqref{eq:subgradients},\eqref{eq:vectors}}{=}\gamma G_x,\notag\\
 q-u&\overset{\eqref{eq:subgradients},\eqref{eq:vectors}}{=}\tau'U,
 &b-B_f&\overset{\eqref{eq:vectors}}{=}s_q-s_u+e.\label{eq:geometry}
\end{align}
The third relation remains valid at $\tau'=0$.
In particular, $p-u=\gamma G_x+\tau'U-\tau G_p$.
Set
\begin{equation}
 h\eqdef\tfrac\gamma{\tau-\gamma},\qquad \theta\eqdef\tfrac\gamma\tau.
 \label{eq:weights}
\end{equation}
Both lie in $(0,1]$, with $\theta\leq\frac12$.
Define the sum of the six weighted residuals:
\begin{align}
 \mathcal I\eqdef{}&(1-\theta)I_f(p,u)+\theta I_f(x,u)+\theta I_f(p,x)\notag\\
 &+(1-h)I_R(p,u;s_u)+hI_R(q,u;s_u)+hI_R(p,q;s_q).
 \label{eq:residuals}
\end{align}
Thus $\mathcal I\overset{\eqref{eq:If},\eqref{eq:Ig}}{\geq}0$. No value of $\RR(x)$ is used.

We first collect the quadratic terms that remain when there is no noise. Define
\begin{align}
 \mathcal Q\eqdef{}&\frac\tau2\norm a^2-\frac{\gamma\tau}{\tau-\gamma}\ip ab
 +\frac{\gamma^2}{\tau-\gamma}\norm b^2
 +\frac{\gamma^2}{\tau-\gamma}\ip{a-b}{B_f}\notag\\
 &+\frac1{2L}\norm{A_f}^2-\frac\gamma{L\tau}\ip{A_f}{B_f}
 +\frac\gamma{L\tau}\norm{B_f}^2.\label{eq:Q}
\end{align}
To verify this form directly, observe that, in the weighted smooth residuals, the function values add up to $\ff(p)-\ff(u)$, with the remaining inner-product part
\[
 \ip{g_u}{p-u}+\theta\ip{g_x-g_u}{p-x}
\]
and squared-gradient part
\[
 \tfrac1{2L}\bigl((1-\theta)\norm{A_f}^2+\theta\norm{B_f}^2
                       +\theta\norm{A_f-B_f}^2\bigr).
\]
In the weighted nonsmooth residuals, the function values add up to $\RR(p)-\RR(u)$, with the remaining inner-product part
\[
 \ip{s_u}{p-u}+h\ip{s_q-s_u}{p-q}.
\]
Add the envelope quadratics $\tau\norm{G_p}^2/2-\tau'\norm U^2/2$ and substitute \eqref{eq:geometry}. Expansion gives
\begin{align}
 \Psi_\tau(x)-\Psi_{\tau'}(q)-\mathcal I
 \overset{\eqref{eq:proxidentity},\eqref{eq:geometry},\eqref{eq:residuals},\eqref{eq:Q}}{=}
 \mathcal Q+\gamma\ip Ue+\frac{\gamma\tau}{\tau-\gamma}\ip ae
                   -\frac{\gamma^2}{\tau-\gamma}\ip be.\label{eq:rawidentity}
\end{align}
The noise terms arise from $s_q-s_u=b-B_f-e$. All $U$-terms without $e$ cancel.
Set
\begin{equation}
 K\eqdef\frac{\gamma^2(\tau-\gamma)}{\tau(\tau-L\gamma^2)},\quad
 C\eqdef\frac{\gamma^2(\tau-2\gamma)}{2(\tau-\gamma)^2},\quad
 J\eqdef\frac{\theta(2-\theta-L\gamma)}{2L},\quad
 \delta\eqdef LK,\quad \eta\eqdef\frac{L\gamma^2}{\tau-L\gamma^2}.\label{eq:squarecoeff}
\end{equation}
We claim the exact identity
\begin{align}
 &\Psi_\tau(x)+\ip{B(x)}e+c\norm e^2-\Psi_{\tau'}(q)\notag\\
 &\quad\overset{\eqref{eq:rawidentity},\eqref{eq:squarecoeff},\eqref{eq:weights},\eqref{eq:Bc}}{=}\mathcal I
 +\frac\tau2\norm{a-h[b-\theta(B_f+e)]}^2
 +\frac1{2L}\norm{A_f-\theta B_f-\delta e}^2\notag\\
 &\qquad\quad+C\norm{b-B_f-e}^2+J\norm{B_f-\eta e}^2.\label{eq:certificate}
\end{align}
Here $c=\frac{\gamma^2}{2(\tau-L\gamma^2)}$ and $B(x)=-\gamma\frac{x-p}{\tau}-K\big(\nabla \ff(p)-\nabla \ff(x)\big)$; when $L\gamma\le1$, these are $c$ in \eqref{eq:c} and $B(x)$ in \eqref{eq:Bc}.
The following two verifications of the coefficients prove this identity.\medskip

\noindent\textbf{Terms independent of $e$.\ }
The coefficients obtained by expanding the four squares in \eqref{eq:certificate} are
\begin{center}\small
\begin{tabular}{lll}\toprule
Term & Expanded coefficient & Coefficient in $\mathcal Q$\\\midrule
$\norm a^2$ & $\tau/2$ & $\tau/2$\\
$\ip ab$ & $-\tau h$ & $-\gamma\tau/(\tau-\gamma)$\\
$\ip a{B_f}$ & $\tau h\theta$ & $\gamma^2/(\tau-\gamma)$\\
$\norm b^2$ & $\tau h^2/2+C$ & $\gamma^2/(\tau-\gamma)$\\
$\ip b{B_f}$ & $-\tau h^2\theta-2C$ & $-\gamma^2/(\tau-\gamma)$\\
$\norm{A_f}^2$ & $1/(2L)$ & $1/(2L)$\\
$\ip{A_f}{B_f}$ & $-\theta/L$ & $-\gamma/(L\tau)$\\
$\norm{B_f}^2$ & $\tau h^2\theta^2/2+\theta^2/(2L)+C+J$ & $\gamma/(L\tau)$\\\bottomrule
\end{tabular}
\end{center}
Each simplification uses the definitions of $h,\theta,C,J$. Thus the noise-free part is precisely $\mathcal Q$.\medskip

\noindent\textbf{Terms involving $e$.\ }
By \eqref{eq:geometry} and \eqref{eq:vectors}, $B(x)=-\gamma(U+a)-K(A_f-B_f)$.
Consequently, since $\frac{\gamma\tau}{\tau-\gamma}-\gamma=\frac{\gamma^2}{\tau-\gamma}=\gamma h$,
\begin{align}
 &\Psi_\tau(x)+\ip{B(x)}e+c\norm e^2-\Psi_{\tau'}(q)-\mathcal I\notag\\
 &\quad\overset{\eqref{eq:rawidentity},\eqref{eq:Bc},\eqref{eq:squarecoeff},\eqref{eq:weights}}{=}
 \mathcal Q+\gamma h\ip{a-b}e-K\ip{A_f-B_f}e+c\norm e^2.\label{eq:noisepolynomial}
\end{align}
The coefficients of $\ip ae$, $\ip be$, $\ip{A_f}e$, and $\ip{B_f}e$
in the proposed squares are, respectively,
\begin{align*}
 \tau h\theta\overset{\eqref{eq:weights}}{=}\gamma h,&\qquad
 -\gamma h^2-2C\overset{\eqref{eq:weights},\eqref{eq:squarecoeff}}{=}-\gamma h,\\
 -\delta/L\overset{\eqref{eq:squarecoeff}}{=}-K,&\qquad
 \gamma h^2\theta+K\theta+2C-2J\eta \overset{\eqref{eq:weights},\eqref{eq:squarecoeff}}{=}K. 
\end{align*}
For example, $\gamma h^2\theta+2C=\gamma^2/\tau$ and
$K\theta-2J\eta=\gamma^3(L\gamma-1)/[\tau(\tau-L\gamma^2)]$, proving the last equality. The squared-noise coefficient is
\begin{align*}
 \frac\tau2h^2\theta^2+C+\frac{\delta^2}{2L}+J\eta^2
 &\overset{\eqref{eq:squarecoeff},\eqref{eq:weights}}{=}
 \frac{\gamma^2}{2\tau}+\frac{L\gamma^4}{2\tau(\tau-L\gamma^2)}
 \overset{\eqref{eq:c}}{=}c.
\end{align*}
These checks and \eqref{eq:noisepolynomial} prove \eqref{eq:certificate}.
They are equalities between the scalar coefficients of polynomials in inner products, and they therefore hold in any real Hilbert space.

Finally, $C,J\geq0$, $\mathcal I\geq0$, and the other square coefficients are positive. Also, $\tau-L\gamma^2>0$:
with $\tau\geq2\gamma$ and $L\gamma\leq2-\gamma/\tau$, it is at least $(\tau-\gamma)^2/\tau>0$. 
Dropping the nonnegative terms in \eqref{eq:certificate} proves \eqref{eq:local} with the above $B(x)$ and $c$, and hence Lemma~\ref{lem:local} when $L\gamma\le1$.
At $\tau'=0$, we have $u=q$ and $s_u=s_q$ by construction, so that all terms remain well defined and the third square has a zero coefficient.

\subsection{The case \texorpdfstring{$L\gamma>1$}{L gamma > 1}}\label{app:long}
Let $1<L\gamma\le2-\theta$. In \eqref{eq:c} and \eqref{eq:Bc}, we then have $\gc=\frac1L$, and
\begin{align}
 &d\overset{\eqref{eq:c}}{=}\tau'+\tfrac1L,\qquad c\overset{\eqref{eq:c}}{=}\frac{L\gamma^2}{2(1+L\tau')},\qquad \mu=\frac{L\gamma-1}{1+L\tau'}>0,\qquad
 K'\eqdef\frac{\gamma(\tau-\gamma)}{\tau(1+L\tau')},\notag\\
 &B(x)\overset{\eqref{eq:Bc}}{=}-\gamma\frac{x-p}{\tau}-K'\big(\nabla \ff(p)-\nabla \ff(x)\big)-\mu\Big(\gamma\frac{x-p}{\tau}-(x-\qc)\Big).\label{eq:longconstants}
\end{align}
We keep the points, vectors, and constants defined above,
and we add the noiseless output $\qc$ of \eqref{eq:localpoints}, with
\begin{equation}
 s_{\qc}\eqdef(x-\qc)/\gamma-g_x\overset{\eqref{eq:localpoints}}{\in}\partial \RR(\qc),\qquad v\eqdef s_q-s_{\qc},\qquad\text{so that}\quad
 \qc-q\overset{\eqref{eq:subgradients}}{=}\gamma(e+v).\label{eq:longpoints}
\end{equation}
The certificate uses three more subgradient inequalities, between $q$ and $\qc$ in both directions and between $\qc$ and $u$: with $\mathcal I$ as in \eqref{eq:residuals}, define
\begin{align}
 \mathcal I'&\eqdef\mathcal I+\mu\,\mathcal I^+,\notag\\
 \mathcal I^+&\eqdef I_R(q,\qc;s_{\qc})+hI_R(\qc,q;s_q)+(1-h)\big[I_R(\qc,u;s_u)-I_R(q,u;s_u)\big].\label{eq:longresiduals}
\end{align}
In $\mathcal I'$, the residual $I_R(q,u;s_u)$ has the weight $h-\mu(1-h)=\frac{1}{1+L\tau'}>0$ and the other eight residuals have nonnegative weights, so that $\mathcal I'\geq0$. The values of $\RR$ cancel in $\mathcal I^+$, and
\begin{equation}
 \mathcal I^+\overset{\eqref{eq:Ig},\eqref{eq:longresiduals}}{=}\ip{s_{\qc}-hs_q-(1-h)s_u}{\qc-q}
 \overset{\eqref{eq:longpoints},\eqref{eq:geometry}}{=}\gamma\ip{(1-h)(b-B_f-e)-v}{e+v}.\label{eq:Iplus}
\end{equation}
Moreover, $\gamma\frac{x-p}{\tau}-(x-\qc)=\gamma(G_p-g_x-s_{\qc})=\gamma(a-b+e+v)$, so that $B(x)=-\gamma(U+a)-K'(A_f-B_f)-\gamma\mu(a-b+e+v)$. Subtracting $\mu\mathcal I^+$ from \eqref{eq:rawidentity} and adding the noise terms therefore gives
\begin{align}
 &\Psi_\tau(x)+\ip{B(x)}e+c\norm e^2-\Psi_{\tau'}(q)-\mathcal I'\notag\\
 &\quad\overset{\eqref{eq:rawidentity},\eqref{eq:Iplus}}{=}
 \mathcal Q+\gamma h\ip{a-b}e-K'\ip{A_f-B_f}e+c\norm e^2\notag\\
 &\qquad\quad-\gamma\mu\ip{a-b+e+v}e-\gamma\mu\ip{(1-h)(b-B_f-e)-v}{e+v}.\label{eq:longpolynomial}
\end{align}
Set
\begin{equation}
 \theta'\eqdef\frac{\tau-L\gamma^2}{\tau(1+L\tau')},\qquad \delta'\eqdef LK',\qquad
 \eta'\eqdef\frac{L\gamma}{1+L\tau'},\qquad N\eqdef\frac{\mu(\tau+L\gamma\tau')}{2(1+L\tau')},\label{eq:longsquarecoeff}
\end{equation}
and keep $C$ and $J$ as in \eqref{eq:squarecoeff}. We claim the exact identity
\begin{align}
 &\Psi_\tau(x)+\ip{B(x)}e+c\norm e^2-\Psi_{\tau'}(q)\notag\\
 &\quad\overset{\eqref{eq:longpolynomial},\eqref{eq:longsquarecoeff},\eqref{eq:squarecoeff}}{=}\mathcal I'
 +\frac\tau2\norm{a-h[b-\theta B_f-\theta'e]}^2
 +\frac1{2L}\norm{A_f-\theta B_f-\delta'e}^2\notag\\
 &\qquad\quad+C\norm{b-B_f-e-\tfrac\mu hv}^2+N\norm v^2+J\norm{B_f-\eta'e}^2.\label{eq:certificatelong}
\end{align}
It is proved by comparing the coefficients of the right-hand side of \eqref{eq:longpolynomial} and of the five squares.
The terms involving neither $e$ nor $v$ are those of the case $L\gamma\le1$: they form $\mathcal Q$ on both sides, by the first table above. The terms involving $v$ and not $e$ are $-\gamma\mu(1-h)\ip{b-B_f}v+\gamma\mu\norm v^2$ in \eqref{eq:longpolynomial} and $-2C\frac\mu h\ip{b-B_f}v+\big(C\frac{\mu^2}{h^2}+N\big)\norm v^2$ in the squares; they agree, since $2C/h=\gamma(1-h)$, $C/h^2=\tau'/2$, and $N=\mu\big(\gamma-\frac{\mu\tau'}{2}\big)$. Finally, the coefficients of $\ip ae$, $\ip be$, $\ip{A_f}e$, $\ip{B_f}e$, $\ip ve$, and $\norm e^2$ in \eqref{eq:longpolynomial},
on the left, and in the squares, on the right, are, respectively,
\begin{align}
 \gamma h-\gamma\mu&\overset{\eqref{eq:weights},\eqref{eq:longconstants},\eqref{eq:longsquarecoeff}}{=}\tau h\theta',\notag\\
 -\gamma h+\gamma\mu h&\overset{\eqref{eq:weights},\eqref{eq:squarecoeff},\eqref{eq:longconstants},\eqref{eq:longsquarecoeff}}{=}-\tau h^2\theta'-2C,\notag\\
 -K'&\overset{\eqref{eq:longsquarecoeff}}{=}-\delta'/L,\notag\\
 K'+\gamma\mu(1-h)&\overset{\eqref{eq:weights},\eqref{eq:squarecoeff},\eqref{eq:longconstants},\eqref{eq:longsquarecoeff}}{=}\tau h^2\theta\theta'+\theta K'+2C-2J\eta',\notag\\
 \gamma\mu(1-h)&\overset{\eqref{eq:weights},\eqref{eq:squarecoeff}}{=}2C\mu/h,\notag\\
 c-\gamma\mu h&\overset{\eqref{eq:weights},\eqref{eq:squarecoeff},\eqref{eq:longconstants},\eqref{eq:longsquarecoeff}}{=}\tfrac\tau2h^2\theta'^2+\tfrac{\delta'^2}{2L}+C+J\eta'^2.\notag
\end{align}
Each of these equalities follows from the definitions \eqref{eq:weights}, \eqref{eq:squarecoeff}, \eqref{eq:longconstants}, and \eqref{eq:longsquarecoeff}; for example, $\gamma h-\gamma\mu=\frac{\gamma(\tau-L\gamma^2)}{(\tau-\gamma)(1+L\tau')}=\tau h\theta'$. These checks and \eqref{eq:longpolynomial} prove \eqref{eq:certificatelong}.

Finally, $C,J,N\geq0$, $\mathcal I'\geq0$, and the other square coefficients are positive.
Dropping the nonnegative terms in \eqref{eq:certificatelong} proves Lemma~\ref{lem:local} when $L\gamma>1$.
At $\tau'=0$, we have $u=q$, $s_u=s_q$, and $h=1$, so that all terms remain well defined and the third square has a zero coefficient. At $L\gamma=1$, we have $\mu=0$, the constants
\eqref{eq:longconstants} and \eqref{eq:longsquarecoeff} take the values of the case $L\gamma\le1$, and \eqref{eq:certificatelong} reduces to \eqref{eq:certificate}.


\section{Proofs of Theorems~\ref{thm:convex}, \ref{thm:scalarconstruction}, and~\ref{thm:secondorder}}\label{app:proofs}

\subsection{Proof of Theorem~\ref{thm:convex}}\label{app:proofconvex}
By \eqref{eq:denominators}, the noise term of index $t$ in \eqref{eq:cost} is $\frac{\sigma^2}{2}\chi(\gamma_t,S_{t+1})$, where, for $z>0$ and $Y\ge0$,
 $\chi(z,Y)\eqdef\frac{z^2}{2Y+\omega(z)}$,
with $\omega$ as in \eqref{eq:omega}. The function $\omega$ is positive and continuously differentiable, with $\omega'\big(\frac1L\big)=0$, so that $\chi$ is continuously differentiable. On the set $\{Lz\le1\}$, $\chi(z,Y)=\phi(r)$ with $r=\frac{z^2}{2(z+Y)}$ and $\phi(r)=\frac{r}{1-Lr}$, as in \eqref{eq:convexcomposition}. On this set, $0<r\leq\frac z2\leq\frac1{2L}$; the function $r$ is a convex quadratic-over-linear function of $(z,Y)$, and $\phi$ is increasing and strictly convex on this interval, so that $\chi$ is convex. On the set $\{Lz\ge1\}$, $\chi(z,Y)=\frac{z^2}{2Y+1/L}$ is a quadratic-over-linear function, which is convex. Hence, along any segment, $\chi$ is a continuously differentiable function of one variable, whose derivative is nondecreasing on each side of the point, if any, where the segment crosses the hyperplane $\{Lz=1\}$; this derivative is therefore nondecreasing, and $\chi$ is convex. Since $(\gamma_t,S_{t+1})$ depends linearly on $\boldsymbol{\gamma}$, and since the initial-error term is convex and smooth, the cost is convex and continuously differentiable on $(0,\infty)^T$.

For strictness, consider a line through a positive schedule, and let $k$ be the largest index of a nonzero component of its direction. Along this line, $S_{k+1}$ is constant and $\gamma_k$ varies at a nonzero rate, and the noise term of index $k$ is a strictly convex function of $\gamma_k$: for a fixed $Y$, the continuously differentiable function $z\mapsto\chi(z,Y)$ is strictly convex on $\{Lz\ge1\}$, where it is a positive quadratic function, and on $\{Lz\le1\}$, where $r$ is a convex function of $z$ with the positive derivative $\frac{z(z+2Y)}{2(z+Y)^2}$ and $\phi$ is increasing and strictly convex. The other terms being convex, the cost is strictly convex along every line, hence strictly convex.

We now prove that a positive minimizer exists. Fix a positive schedule, let $\bar F$ be its cost, and consider the positive schedules whose cost is at most $\bar F$. Their total stepsize is at least $\frac{D^2}{4\bar F}$, because of the initial-error term, and their steps are bounded: since $\omega\le\frac1L$, the noise term of index $t$ is at least $\frac{\sigma^2\gamma_t^2}{2(2S_{t+1}+1/L)}$, so that $\gamma_t^2\le\frac{2\bar F}{\sigma^2}\big(2S_{t+1}+\frac1L\big)$, which bounds $\gamma_{T-1}$, then $\gamma_{T-2}$, and so on backward. Extend the cost to the schedules with nonnegative steps and a positive total stepsize by assigning a zero noise term to a zero step. This extension is continuous, because, for $L\gamma_i\le1$,
\[
 0\leq\frac{\gamma_i^2}{d_i}
 \overset{\eqref{eq:denominators}}{\leq}\frac{\gamma_i^2}{\omega(\gamma_i)}=\frac{\gamma_i}{2-L\gamma_i}\leq\gamma_i.
\]
By compactness, the extended cost has a minimizer over the schedules with nonnegative steps, and its total stepsize is positive. No coordinate with a positive tail can be zero there:
the right derivative of the cost with respect to such a zero coordinate $\gamma_k$ is
$-D^2/(4S_0^2)-\sigma^2\sum_{i<k}\gamma_i^2/d_i^2<0$.

There can be no trailing zero block either. If $\gamma_m=a>0$ is the last positive coordinate and $m<T-1$, replace $(\gamma_m,\gamma_{m+1})$ by $(a-\varepsilon,\varepsilon)$. All earlier sums and the initial-error term stay fixed, and the two affected noise fractions become $\frac{(a-\varepsilon)^2}{2\varepsilon+\omega(a-\varepsilon)}+\frac{\varepsilon}{2-L\varepsilon}$. Their derivative at $\varepsilon=0+$ is
$-\frac4{(2-La)^2}+\frac12$ if $La\le1$ and $-2La(1+La)+\frac12$ if $La\ge1$; it is negative in both cases, which contradicts optimality. The minimizer is therefore positive, and strict convexity proves uniqueness.

Differentiating \eqref{eq:cost} gives the partial derivatives \eqref{eq:gradient}: the denominator $d_i$ depends on $\gamma_k$ through $S_{i+1}$ if $i<k$ and through $\omega(\gamma_k)$ if $i=k$, and $2d_k-\gamma_k\omega'(\gamma_k)=2\big(2S_{k+1}+\min\{\gamma_k,\frac1L\}\big)$, whether $L\gamma_k\le1$ or $L\gamma_k\ge1$. Since the cost is convex and differentiable on the open set $(0,\infty)^T$, a positive schedule is its minimizer if and only if its gradient vanishes there, which is \eqref{eq:kkt}.

\subsection{Proof of Theorem~\ref{thm:scalarconstruction}}\label{app:proofscalar}
For arguments at most $\frac1L$, the functions \eqref{eq:scalarfunctions} are $P_Y(a)=\frac{2a(a+Y)}{[2(a+Y)-La^2]^2}$ and $B_Y(b)=\frac{b(2Y-b)}{(2Y-Lb^2)^2}$; for arguments at least $\frac1L$, they are $P_Y(a)=\frac{a(a+2Y+1/L)}{(2Y+1/L)^2}$ and $B_Y(b)=\frac{b}{2Y-2b+1/L}$; they are continuous. For $Y>0$, differentiation gives
\begin{equation}
 P_Y'(a)\overset{\eqref{eq:scalarfunctions}}{=}
 \frac{4Y(a+Y)+2La^2(2a+3Y)}{[2(a+Y)-La^2]^3}>0\ \ \text{if }La<1,\qquad
 P_Y'(a)\overset{\eqref{eq:scalarfunctions}}{=}\frac{2a+2Y+\frac1L}{(2Y+\frac1L)^2}>0\ \ \text{if }La>1,
 \label{eq:scalarmonotonicity}
\end{equation}
and $P_Y(a)\to+\infty$ as $a\to+\infty$.
Also, $P_Y(b)<B_Y(b)$ by \eqref{eq:tiederivative} and \eqref{eq:tiederivativelong},
with tail $Y-b\geq0$. Thus the equation $P_Y(a)=B_Y(b)$ has a unique solution $a>b$, so that the construction is well defined and the schedule $\Gamma(s)$ is strictly decreasing.
By \eqref{eq:adjacentstationarity}, applied to the schedule $\Gamma(s)$,
all the components of the gradient of the cost at $\Gamma(s)$ are equal to $\Lambda(s)$.
Hence, whenever $\Lambda(s)=0$, the conditions \eqref{eq:kkt} hold and $\Gamma(s)$ is the minimizer.

It remains to justify the scalar search.
The monotonicity in \eqref{eq:scalarmonotonicity} implies that $\Gamma(s)$, and hence $\Lambda(s)$, are continuous. Moreover, $\Gamma(s)\to0$ as $s\downarrow0$.
To prove the latter backward, suppose $b,Y\to0$. Since $b\leq Y$,
$B_Y(b)\overset{\eqref{eq:scalarfunctions}}{\leq}\tfrac1{(2-Lb)^2}\longrightarrow\tfrac14$.
For any fixed $a>0$, however, $P_Y(a)\to a\big(a+\min\{a,\tfrac1L\}\big)/\omega(a)^2\geq1/2$. No positive subsequential limit of the preceding step is therefore compatible with the construction.
This completes the induction.
Formula \eqref{eq:gradient} now gives, for $s\le\frac1L$, $\Lambda(s)\leq\sigma^2-D^2/(4S_0(s)^2)\to-\infty$.

For $s\ge\frac1L$, all the steps of $\Gamma(s)$ are at least $\frac1L$; writing $S_0$ and $S_1$ for its tails,
\begin{equation*}
 \Lambda(s)=\partial_0F(\Gamma(s))\overset{\eqref{eq:gradient}}{=}-\frac{D^2}{4S_0^2}+\frac{\sigma^2\Gamma_0}{2S_1+\frac1L}
 \geq\frac{1}{S_0^2}\Big(\sigma^2\Gamma_0\big(2\Gamma_0-\tfrac1L\big)-\frac{D^2}{4}\Big),
\end{equation*}
because, as in \eqref{eq:firststationarity}, the minimum of $\frac{S_0^2}{2S_1+1/L}=\frac{(\Gamma_0+S_1)^2}{2S_1+1/L}$ over $S_1\ge0$ is $2\Gamma_0-\frac1L$, attained at $S_1=\Gamma_0-\frac1L$.
Since $\Gamma_0\ge s$, we have $\Lambda(s)>0$ for $s$ large enough, and continuity gives a zero.
Two distinct zeros would yield two distinct schedules satisfying \eqref{eq:kkt}, contradicting the uniqueness of the minimizer.
Thus the zero is unique, and $\Lambda$ is negative below it and positive above it.
This validates bisection and completes the proof.

\subsection{Proof of Theorem~\ref{thm:secondorder}}\label{app:proofsecondorder}
The lower comparison is \eqref{eq:globalorder}. We prove the upper bound at $\widehat{\boldsymbol{\gamma}}$.

For $T=1$, its largest step is $D/\sigma$.
For $T\geq2$, the monotonicity stated in Theorem~\ref{thm:ppaopt} and $Q_{T-1}\leq1$ give
\begin{equation*}
 \max_t\widehat\gamma_t
 \overset{\eqref{eq:ppascale},\eqref{eq:ppaschedule},\eqref{eq:qrec}}{=}
 \frac{D\sqrt{Q_T}}{2\sigma(1-Q_{T-1}/4)}
 \leq\frac{2D\sqrt{Q_T}}{3\sigma}
 \leq\frac{D\sqrt{Q_T}}{\sigma}.
\end{equation*}
The hypothesis of the theorem therefore ensures that all the steps are at most $\frac1L$, also when $T=1$.

Let $\widehat S_t$ denote the tails of $\widehat{\boldsymbol{\gamma}}$.
For $n=T-t\geq2$,
\begin{equation}
 \frac{\widehat S_{t+1}}{\widehat S_t}
 \overset{\eqref{eq:ppaschedule}}{=}1-\frac{Q_{n-1}}2
 \leq\left(1-\frac{Q_{n-1}}4\right)^2
 \overset{\eqref{eq:qrec}}{=}\left(\frac{Q_n}{Q_{n-1}}\right)^2.
 \label{eq:ppatailratio}
\end{equation}
Multiplying these inequalities yields
\begin{equation}
 \widehat S_t\overset{\eqref{eq:ppatailratio}}{\leq}
 \widehat S_0\left(\frac{Q_T}{Q_{T-t}}\right)^2.
 \label{eq:ppatailbound}
\end{equation}
Put $\widehat r_t\eqdef\widehat\gamma_t^2/(2\widehat S_t)$. At the last stage, \eqref{eq:ppatailbound} gives $\widehat r_{T-1}\leq\widehat S_0Q_T^2/2$. At any earlier stage, \eqref{eq:ppaschedule} and \eqref{eq:ppatailbound} give the smaller upper bound $\widehat S_0Q_T^2/[8(1-Q_{n-1}/4)^2]\leq2\widehat S_0Q_T^2/9$.
Also, the optimal normalized noise cost is $Q_T$. Hence
\begin{align}
 \max_t\widehat r_t&\overset{\eqref{eq:ppascale},\eqref{eq:ppatailbound},\eqref{eq:ppaschedule}}{\leq}
  \frac{DQ_T^{3/2}}{2\sigma},\qquad
 \sum_t\widehat r_t\overset{\eqref{eq:shape},\eqref{eq:bellman},\eqref{eq:ppascale}}{=}
 \frac{D\sqrt{Q_T}}{2\sigma}.\label{eq:rsum}
\end{align}
Therefore
\begin{equation}
 \sum_t\widehat r_t^2\overset{\eqref{eq:rsum}}{\leq}\frac{D^2Q_T^2}{4\sigma^2}.
 \label{eq:rsquares}
\end{equation}

Equations \eqref{eq:transferparameter}, 
\eqref{eq:rsum},
and \eqref{eq:rsquares} bound the curvature penalty \eqref{eq:generalpenalty} of $\widehat{\boldsymbol{\gamma}}$ by the last term of \eqref{eq:secondordersandwich}.
Theorem~\ref{thm:exact} and the definitions of the infima then prove \eqref{eq:secondordersandwich}.

Finally, for $y>0$ and $z\geq0$, we have $(y+z)^{-2}\geq y^{-2}-2zy^{-3}$.
Apply this with $y=V_{0,T}$ and $z$ equal to the last term of \eqref{eq:secondordersandwich};
using \eqref{eq:ppavalue} and $b_T=4/Q_T$ yields \eqref{eq:effectivehorizon}. Since $Q_T\to0$, \eqref{eq:basymptotic} then gives \eqref{eq:universalconstant}.

\section*{Acknowledgments}
This work was supported by funding from King Abdullah University of Science and Technology (KAUST):
i) KAUST Baseline Research Scheme, and 
ii) Competitive Research Grant (CRG) Program (award no.\ 6460).

The authors used the AI systems Claude Fable (Anthropic) and GPT-6 Astra (OpenAI) in the derivation and the verification of some of the results and in the writing of the paper. They have checked the results and take full responsibility for the content.

\bibliographystyle{abbrvnat}
\bibliography{biblio}

\end{document}